\documentclass[11pt]{amsart}

\usepackage[leqno]{amsmath}
\usepackage{amsthm}
\usepackage{amsfonts}
\usepackage{amssymb}
\usepackage{eucal}
\usepackage[all]{xy}
\CompileMatrices

\DeclareFontFamily{OT1}{rsfs}{}
\DeclareFontShape{OT1}{rsfs}{n}{it}{<-> rsfs10}{}
\DeclareMathAlphabet{\mathscr}{OT1}{rsfs}{n}{it}

\newcommand{\ie}{{\it i.e.}}

\newcommand{\lead}{\operatorname {lead}}

\newcommand{\disc}{\operatorname {disc}}

\newcommand{\FF}{{\mathbf F}}
\newcommand{\kap}{\kappa}

\newcommand{\RR}{{\mathbf R}}
\newcommand{\ZZ}{{\mathbf Z}}
\newcommand{\PP}{{\mathbf P}}
\newcommand{\QQ}{{\mathbf Q}}
\newcommand{\OO}{{\mathscr O}}

\newcommand{\Spec}{\operatorname {Spec}}

\newcommand{\Gal}{\operatorname {Gal}}

\newcommand{\ord}{\operatorname {ord}}

\newcommand{\chr}{\operatorname {char}}
\newcommand{\K}{K}
\newcommand{\KK}{{\mathcal K}}
\newcommand{\avg}{{\rm{Avg}}}

\DeclareFontEncoding{OT2}{}{} 
  \newcommand{\textcyr}[1]{%
    {\fontencoding{OT2}\fontfamily{wncyr}\fontseries{m}\fontshape{n}%
     \selectfont #1}}
\newcommand{\Sha}{{\mbox{\textcyr{Sh}}}}

\theoremstyle{plain}
\newtheorem{theorem}{Theorem}[section]
\newtheorem{lemma}[theorem]{Lemma}
\newtheorem{corollary}[theorem]{Corollary}

\theoremstyle{definition}
\newtheorem{remark}[theorem]{Remark}
\newtheorem{definition}[theorem]{Definition}
\newtheorem{example}[theorem]{Example}

\newcounter{bean}

\numberwithin{equation}{section}

\newcommand{\eqdef}{:=}

\newcommand{\Atop}[2]{\genfrac{}{}{0pt}{}{#1}{#2}}

\newcommand{\Nm}{\operatorname N}

\newcommand*{\invlim}{\varprojlim}

\title{Root numbers and ranks in positive characteristic}
\author{B. Conrad}
\address{Department of Mathematics,
Univ. of Michigan, Ann Arbor, MI 48109-1109, USA}
\email{bdconrad@umich.edu}

\author{K. Conrad}
\address{Department of Mathematics,
Univ. of Connecticut, Storrs, CT 06269-3009, USA}
\email{kconrad@math.uconn.edu}

\author{H. Helfgott}
\address{D\'epartement de Math\'ematiques et Statistique, 
Universit\'e de Montr\'eal, Montr\'eal, QC H3C 3J7, Canada}
\email{helfgott@dms.umontreal.ca}

\thanks{We thank N. Boston, I. Dolgachev, A.J. de Jong, B. Mazur, 
D. Pollack, B. Poonen,
K. Rubin, A. Silberstein, E. Spiegel, R. Vakil, 
A. Venkatesh, and S. Wong for helpful conversations. 
B.C. thanks the NSF and the Sloan Foundation
for support.}

\begin{document}

\begin{abstract}
For a global field $K$ and an elliptic curve
$\mathscr{E}_{\eta}$ over $K(T)$, Silverman's specialization theorem implies
$
{\rm{rank}}(\mathscr{E}_{\eta}(K(T))) \le 
{\rm{rank}}(\mathscr{E}_t(K))
$
for all but finitely many 
$t \in \PP^1(K)$.
If this inequality is strict
for all but finitely many $t$, 
the elliptic curve $\mathscr{E}_{\eta}$ 
is said to have {\em elevated rank}. 
All known examples of elevated rank for $K = \QQ$
rest on the parity conjecture for elliptic curves over 
$\QQ$, and the examples are all isotrivial.  

Some additional standard conjectures over
$\QQ$ imply that there does not exist 
a non-isotrivial elliptic curve over $\QQ(T)$ with
elevated rank.
In positive characteristic, an analogue of one of 
these additional conjectures is false. 
Inspired by this, for the rational function field
$K = \kappa(u)$ over any finite field
$\kappa$ with characteristic $\ne 2$,  
we construct an explicit 2-parameter family
$E_{c,d}$ of non-isotrivial 
elliptic curves 
over $K(T)$ (depending on arbitrary $c, d \in \kappa^{\times}$)
such that, under the parity conjecture, each $E_{c,d}$ has elevated rank. 
\end{abstract}

\maketitle

\begin{center}
{\it To Mike Artin on his 70th birthday}
\end{center}

\section{Introduction}\label{intro}

Let $K$ be a global field and let 
$\mathscr{E}_{\eta}$ be an elliptic curve over $K(T)$.
This curve uniquely extends to a minimal
regular proper elliptic fibration
$\mathscr{E} \rightarrow \PP^1_\K$.
The group $\mathscr{E}_{\eta}(K(T))$ 
is finitely generated, by the Lang--N\'eron theorem
\cite[Thm.~1]{langneron}. 
(See \cite[\S6]{bnotes} for a proof of the Lang--N\'eron
theorem using the language of schemes.)
For all but finitely many $t \in \PP^1(K)$, 
the specialization $\mathscr{E}_t$ 
of $\mathscr{E}$ at $T = t$ is an elliptic curve over $K$.  
This paper is concerned with a comparison between 
the ranks of $\mathscr{E}_{\eta}(K(T))$ and $\mathscr{E}_t(K)$ 
as $t$ varies.  

By Silverman's specialization theorem 
\cite[Thm.~C]{silrank}, the specialization map
$$
\mathscr{E}_{\eta}(K(T)) = \mathscr{E}(\PP^1_K)
\rightarrow \mathscr{E}_t(K)
$$
at $T = t$ is injective for all but finitely many 
$t \in \PP^1(K)$.
(To be precise, Silverman's theorem only 
applies to non-constant 
$\mathscr{E}_{\eta}$.  Injectivity of the 
specialization map for constant $\mathscr{E}_{\eta}$ is
elementary.)  
Thus, the {\em generic rank} 
$r(\mathscr{E}_{\eta}) 
:= {\rm{rank}}(\mathscr{E}_{\eta}(K(T)))$ 
satisfies 
\begin{equation}\label{ineq}
r(\mathscr{E}_{\eta}) \le 
{\rm{rank}}(\mathscr{E}_t(\K))
\end{equation}
for all but finitely many $t$. 
The elliptic curve ${\mathscr{E}}_{\eta}$
(or the fibration $\mathscr{E} \rightarrow 
\PP^1_K$, or the family $\{\mathscr{E}_t\}_{t \in \PP^1(K)}$)
is said to have 
{\it elevated rank} if (\ref{ineq}) is a strict inequality 
for all but finitely many $t \in \PP^1(K)$.

How are examples of elevated rank constructed? 
The only known technique depends on 
the {\it parity conjecture}: 
for every elliptic curve $E$ over the global field $K$, 
$$
(-1)^{{\rm rank}(E(K))} \stackrel{?}{=} W(E), 
$$
where $W(E)$ is the global root number of $E$. 
The spirit of the parity 
conjecture is that $W(E)$ is supposed to be the sign 
in the 
functional equation of the $L$-function of $E$, but 
such a functional equation 
is not yet known to exist in general.  
Therefore, we adopt the convention that the global root number is 
defined to be the product of local root numbers. 
The local root numbers are defined in all cases  
via representation theory \cite{deligne} and are equal to 1 at
non-archimedean places of good reduction. 
Some convenient formulas for local root numbers 
at non-archimedean places will be recalled in Theorem \ref{localrtnum}.  
The analytic and representation-theoretic descriptions of $W(E)$  
are known to agree when $K$ is $\QQ$ or a global function field, by 
work of Deligne, Drinfeld, Wiles, and others.  In particular, since our 
focus in this paper will be the cases when $K = \QQ$ or when
$K$ is a rational function field over a finite field, 
the reader can think about $W(E)$ in either way. 

To find elliptic curves with elevated rank, one 
tries to construct $\mathscr{E} \rightarrow \PP^1_K$ such that 
$W(\mathscr{E}_t)$ has opposite sign 
to $(-1)^{r(\mathscr{E}_{\eta})}$ with at most finitely many exceptions.  
That is, we want 
\begin{equation}\label{Wnot}
W(\mathscr{E}_t) = -(-1)^{r(\mathscr{E}_{\eta})}
\end{equation}
for all but finitely many $t \in \PP^1(K)$.  
Assuming the parity conjecture for elliptic curves over $K$, 
(\ref{ineq}) and (\ref{Wnot}) 
imply that (\ref{ineq}) is a strict 
inequality for all but finitely many $t$, so 
the $\mathscr{E}_t$'s have elevated rank.

Because of the role of the parity conjecture in this strategy,
all known examples of elevated rank are, strictly speaking, 
conditional.  Moreover, this idea has so far only 
been carried out when $K = \QQ$. 
The first (conditional) examples of elevated rank 
were found by Cassels and 
Schinzel \cite{casselsschinzel}.  These are quadratic twists over 
$\QQ(T)$ of the elliptic curve $y^2 = x^3 - x$: 
\begin{equation}\label{csex}
\mathscr{E}_{n,\eta} : n(1+T^4)y^2 = x^3 - x,
\end{equation}
where $n$ is a (squarefree) positive integer satisfying
$n \equiv 7 \bmod 8$.
Each ${\mathscr{E}}_{n,\eta}$ should have elevated rank, because the 
group $\mathscr{E}_{n,\eta}(\QQ(T))$ 
has rank 0 and $W(\mathscr{E}_{n,t}) = -1$ 
for every $t \in \QQ$.
We exclude $t = \infty$ because $\mathscr{E}_{n,\infty}$ is not
smooth.
(Although $\mathscr{E}_{n,t}(\QQ)$ 
should have a point of infinite order for each $t \in \QQ$,  
there cannot be an algebraic formula
for a point of infinite 
order on $\mathscr{E}_{n,t}(\QQ)$
as $t$ varies through an infinite subset
of $\QQ$, since such a formula would give an element of 
infinite order in the group $\mathscr{E}_{n,\eta}(\QQ(T))$.) 
More generally, for any 
elliptic curve $E_{/\QQ}$, Rohrlich
\cite[Prop. 9]{rohrlich1} proved that there is a quadratic twist 
$\mathscr{E}_{\eta}$ of $E_{/\QQ(T)}$ 
by a quartic irreducible in $\QQ[T]$
such that $\mathscr{E}_{\eta}(\QQ(T))$ has rank 0 
and $W(\mathscr{E}_t) = -1$ for every $t \in \QQ$.

Nekov\'a\v{r} has proved 
the parity conjecture for any elliptic curve over $\QQ$ 
with finite Tate--Shafarevich group 
\cite{nekovar}, \cite[p.~463]{rubsil}, 
but this does not make any examples of 
elevated rank over $\QQ$ unconditional, 
since there are no non-constant families 
$\mathscr{E} \rightarrow \PP^1_{\QQ}$
such that  $\mathscr{E}_t$ is known to 
have a finite Tate--Shafarevich group for 
all but finitely many $t \in \PP^1(\QQ)$.   
Similarly, the recent work of 
Kato and Trihan \cite{katotrihan} 
(as well as earlier work of Artin--Tate, Milne,
Schneider, and others) on the Birch and Swinnerton-Dyer 
conjecture in characteristic $p$ 
does not have any impact
on the conditional nature of the parity conjecture in 
characteristic $p$ as it is applied
to the examples considered in this paper.

The examples of Cassels--Schinzel and Rohrlich
over $\QQ(T)$ are quadratic twists. 
The appeal of quadratic twists
is that there are simple formulas 
that describe the variation of root numbers under 
quadratic twists over $\QQ$ \cite[Cor. to Prop. 10]{rohrlich3}, 
\cite[Thm.~7.2]{rubsil}.  However, a family of quadratic 
twists exhibits no ``geometric'' variation:  it is isotrivial 
(that is, $j(\mathscr{E}_{\eta}) \in K$, with no $T$-dependence), and 
conversely any isotrivial family is either a family of quadratic 
twists or is a family of quartic (resp. cubic or
sextic) twists 
with $j(\mathscr{E}_{\eta}) = 1728$ (resp. $j(\mathscr{E}_{\eta}) = 0$).

The main question we address in this paper is the following: 
for a global field $K$, does there exist
a {\em non-isotrivial} elliptic curve
over $K(T)$ with elevated rank?
In Appendix \ref{avgrtnum}, 
we explain why some standard conjectures
over $\QQ$ imply that the answer to this question for $K = \QQ$
is {\it no}.
There are natural analogues of these standard conjectures
over a rational function field
$\kappa(u)$ over a finite field $\kappa$, 
but (as we explain in Appendix \ref{contrsec}) one of these 
conjectures is false over $\kappa(u)$. 
This suggests that our question might have an affirmative 
answer in the function field case.

Here is our example.  Let $\kappa$ be a finite field
with characteristic $p \ne 2$, and choose any  
$c, d \in \kappa^\times$.  
Let $F = \kappa(u)$ and consider the elliptic curve
\begin{equation}\label{mainfamily}
\mathscr{E}_{\eta} : 
y^2 = x^3 + (c(T^2+u)^{2p}+du)x^2 -(c(T^2+u)^{2p}+du)^3x
\end{equation}
over $F(T)$.
The Weierstrass model 
(\ref{mainfamily}) over $F(T)$ has the form $y^2 = x^3 + Ax^2 - A^3x$. 
The $j$-invariant $j(\mathscr{E}_{\eta}) \in F(T)$ is not in $F$, so
$\mathscr{E}_{\eta/F(T)}$ is non-isotrivial.
An inspection of the poles of $j(\mathscr{E}_{\eta})$ on
$\PP^1_F$ shows that changing $(c,d)$ changes
$j(\mathscr{E}_{\eta})$.

Let $\mathscr{E} \rightarrow \PP^1_F$ be the minimal 
regular proper elliptic fibration
with generic fiber $\mathscr{E}_{\eta}$.  For all $t \in \PP^1(F)$, 
the specialization
$\mathscr{E}_t$ of $\mathscr{E}$ at $T = t$ is an elliptic curve.

\begin{theorem}\label{mainthm}
Let $F = \kappa(u)$,
where ${\rm{char}}(\kappa) \ne 2$,
 and fix a choice of $c, d \in \kappa^\times$. 
Let $\mathscr{E}_{\eta/F(T)}$ be as in $(\ref{mainfamily})$, 
depending on 
the choice of $c$ and $d$.
For every $t \in \PP^1(F)$, 
we have $W(\mathscr{E}_t) = 1$. 
If $t \ne \infty$, then 
$\mathscr{E}_t(F)$ has
positive rank.  Moreover, 
${\rm{rank}}(\mathscr{E}_{\eta}(F(T))) = 1$ and 
$$W(\mathscr{E}_t) = -(-1)^{{\rm{rank}}(\mathscr{E}_{\eta}(F(T)))}$$
for all $t \in \PP^1(F)$. 

Thus, if the parity conjecture is true for elliptic curves 
over 
$F$, the $\mathscr{E}_t$'s are a non-isotrivial family with elevated rank. 
\end{theorem}

\begin{remark}\label{infinity}
When $t = \infty$, the fiber ${\mathscr E}_t$ in Theorem \ref{mainthm} 
is the constant elliptic curve 
$y^2 = x^3 - c^3 x$ over $F$.  The elliptic curve
$\mathscr{E}_{\infty}$ therefore 
has global root number 1, and it must have
rank 0 (since $F$ is a function field of genus 0 over
the finite field $\kappa$). 

We expect that
$\mathscr{E}_t(F)$ has rank 2 except for a set of $t \in 
\PP^1(F)$ with 
density 0 (as measured by height), but 
we have no idea how to prove this expectation.
\end{remark}

\begin{remark} 
We did not search for non-isotrivial examples of
elevated rank with generic rank 0 or in characteristic 2, 
but we expect that such examples exist.
(The curve defined by (\ref{mainfamily}) in characteristic 2 
is not smooth.)  Our example has 
$j(\mathscr{E}_{\eta}) \in F(T^p)$, and we expect 
any example of elevated rank in 
characteristic $p$ will have the same property.  
\end{remark}

\begin{remark}
For fixed $p \ne 2$, consider
the algebraic family
$\mathscr{E}_t$ where $\kappa$ and $(c,d)$ vary (in characteristic 
$p$) but
the logarithmic height of $t \in \kappa(u)^{\times}$ 
(\ie, the maximum of the degrees of its
numerator and denominator)
is bounded by some integer $B > 0$.  This is a family parameterized
 by the $\kappa$-points $(c,d,t)$ of
a smooth $\FF_p$-scheme that is determined by
$B$.

For fixed $c, d \in \kappa^{\times}$, the 
assertions in Theorem \ref{mainthm} concerning the fibral root numbers
and the generic rank for the associated
elliptic curve in (\ref{mainfamily}) are unaffected
by replacing $\kappa$ with a finite extension.
(This is also crucial in the proof of 
Theorem \ref{mainthm}; see (\ref{specdig}) and 
the surrounding text.) 
Thus, granting the parity conjecture, 
Theorem \ref{mainthm} implies that there is
a systematic ``rank
gap'' $\ge 1$ between generic and special Mordell--Weil ranks
over the connected components of the family of 
planar Weierstrass models $\mathscr{E}_t$ as $(c,d,t)$ and $\kappa$ vary
with ${\rm{height}}(t) \le B$.
This is a
contrast with a theorem in
\cite[\S9]{deJongkatz} asserting that there is an average ``rank gap''
$\le 1/2$ (or exactly $1/2$
under conjectures of Tate)
between generic and special Mordell--Weil ranks 
of Jacobians of 
certain {\em universal} families of pencils 
of smooth plane curves in characteristic $p$. 
(The pencils considered in
\cite[\S9]{deJongkatz} are induced by 
smooth surfaces in $\PP^1 \times \PP^2$, but the 
closure of (\ref{mainfamily}) in $\PP^1_F \times_F \PP^2_F$ is
not $F$-smooth.) 
\end{remark}

Here is an overview of how we prove Theorem \ref{mainthm}.
To compute the rank of 
a Mordell--Weil group, we wish to use a 2-descent, 
and this is simplest when there is
a rational point of order 2 and we are not in characteristic 2.  
Weierstrass models for such elliptic curves
can always be brought to the form
$y^2 = x^3 + Ax^2 + Bx$ with the 2-torsion point equal to $(0,0)$. 

\underline{Step 1}:  For an odd prime $p$,
consider the non-isotrivial elliptic curve over 
$\FF_p(T)$ given by the Weierstrass model with $A = T$, $B = -T^3$:
$$
E_T : y^2 = x^3 + Tx^2 - T^3x.
$$
For every $t \in \kappa(u)$ with $t \ne 0, -1/4$, the
specialization 
$E_t$ may be considered as an elliptic curve over $\kappa(u)$. 
We will compute 
the reduction type of $E_t$ at every place of $\kappa(u)$.  
When $\chr(\kappa) = 3$, the 2-torsion point (0,0) will 
prevent the intervention of wild ramification.

\underline{Step 2}: We show that $Q_T = (-T,T^2) \in 
E_T(\FF_p(T))$ has infinite order, so $Q_t$ has infinite 
order in $E_t(\kappa(u))$ for every $t \in \kappa(u)$ such that $t \not\in 
\kappa$.  This uses an extension to characteristic $p$ of 
the classical Nagell--Lutz criterion in characteristic 0. 

\underline{Step 3}: Letting 
$h(T) = cT^{2p} + du$, where 
$c, d \in \kappa^\times$, 
we use algebraic properties
of the defining Weierstrass model for
$E_T$ to 
find a simple formula (\ref{WW}) for $W(E_{h(t)})$ 
for {\em every}
$t \in \kappa(u)$.  (Note $h(t) \not\in \kappa$ for 
every $t$.) The formula
for $W(E_{h(t)})$ implies that $W(E_{h(t^2 + u)}) = 1$ for
all $t \in \kappa(u)$; the elliptic curve 
$E_{h(T^2+u)}$ 
is (\ref{mainfamily}).
Our specific choice of $h(T)$ is partly motivated by
the study of the characteristic-$p$ M\"obius function in 
\cite{ccg}.

\underline{Step 4}: As just noted, $\mathscr{E}_{\eta}$ in (\ref{mainfamily}) is 
$E_{h(T^2+u)}$. The point $Q_{h(T^2+u)}$ on this curve has infinite order, and we use 
this point to show that $\mathscr{E}_t(F)$ has positive rank
for all $t \in \PP^1(F) - \{\infty\}$.
A  mixture of
geometric, arithmetic, and cohomological arguments 
is used to prove that the
rank of $\mathscr{E}_{\eta}(F(T))$ is $\le 1$
(so the rank is exactly 1).
The essential inputs are 
the Lang--N\'eron theorem over an algebraic closure 
$\overline{\kappa}$,
the geometry of the locus of bad reduction for
$\mathscr{E}_{\eta}$ over $\PP^1 \times \PP^1$, 
and some arithmetic considerations with the Chebotarev density theorem.
Standard geometric upper bounds on the $F(T)$-rank 
give very large bounds when applied to $\mathscr{E}_{\eta}$.
It therefore 
seems hopeless to calculate the rank of $\mathscr{E}_{\eta}(F(T))$ 
via purely geometric methods, even though the generic-rank conclusion
in Theorem \ref{mainthm} holds over $\overline{\kappa}$. 

Steps 1 and 2 are carried out in \S\ref{sec3},
Step 3 is carried out in \S\ref{rtlocnum},
and Step 4 is carried out in \S\ref{rank0}--\S\ref{selmer}.
The bulk of the work is in Step 4, which the 
geometrically-inclined reader may prefer to read 
directly after Step 1. 
We note that Steps 2 and 3 are logically independent, as are 
Steps 3 and 4.  (Clearly Step 2 is used in Step 4.)
In \S\ref{nagsec}, we
discuss a rank conjecture of Nagao in the context of (\ref{mainfamily}).

In Appendix \ref{avgrtnum}, we review previous work on 
variation of root numbers in
families over $\QQ$.  
The reason to expect the possibility of different behavior in 
positive 
characteristic is explained in Appendix \ref{contrsec}. 
These appendices are
expository, but they should help 
the reader to have the proper perspective on 
our work. 

Our notation is standard, with two exceptions:
$F$ denotes the rational function field 
$\kappa(u)$, where $\kappa$ is a finite field
that is assumed to have
characteristic $\not= 2$ unless otherwise stated,
and in some calculations in a field we shall use the shorthand
$x \sim y$ 
to denote the relation $x = y z^2$ for a
non-zero $z$ (see Definition \ref{funnydef}).

\section{Reduction type and rank for $y^2 = x^3 + Tx^2 - T^3x$}\label{sec3}

We begin with two elementary lemmas concerning reduction types for 
an elliptic curve over the fraction field of a
discrete valuation ring.  We write $\KK$ for the fraction 
field and $v$ for the normalized (\ie, $\ZZ$-valued) 
discrete valuation on $\KK$, with valuation ring $\OO_{\KK}$ and 
residue field $k$.
Both lemmas are standard when
${\rm{char}}(k) \ne 2, 3$, so
a key point of the proofs is to include
the case $\chr(k) = 3$.
For an elliptic curve $E$ over $\KK$, 
we let $\Delta$, $c_4$, and $c_6$ denote the usual parameters
associated to a Weierstrass model for $E$ 
over $\KK$.
(As is well-known, $\Delta \bmod (\KK^{\times})^{12}$, $c_4 \bmod
(\KK^{\times})^4$, and $c_6 \bmod (\KK^{\times})^6$ are
independent of the choice of Weierstrass model.)

\begin{lemma}\label{good}
Let $E$ be an elliptic curve over $\KK$ with 
potentially good reduction.  If there is good reduction, then 
$v(\Delta) \equiv 0 \bmod 12$.  The converse holds in
either of the two following situations:
\begin{enumerate}
\item $\chr(k) \ne 2, 3$,
\item $\chr(k) \ne 2$ and $E(\KK)[2] \ne O$.
\end{enumerate}
\end{lemma}

\begin{proof}
The necessity of
the congruence condition for good reduction is obvious.
When $\chr(k) \ne 2, 3$, all integral Weierstrass
models can be put in the form
$y^2 = x^3 + \alpha x + \beta$, so the sufficiency 
of the congruence condition for good reduction 
in case (1) is proved by direct calculation
using the standard formulas for $j$ and $\Delta$ (in terms of $\alpha$
and $\beta$) and
using the coordinate changes
$(x,y) \mapsto (\gamma^2 x, \gamma^3 y)$
with $\gamma \in \KK^{\times}$.

For sufficiency in case (2), we may suppose $\OO_{\KK}$ is
strictly henselian.  
Thus, $k$ is separably closed with $\chr(k) \ne 2$,
so $\Delta$ must be a square in $\KK^{\times}$ since 
$v(\Delta)$ is even.  By the 2-torsion hypothesis, 
for any Weierstrass $\KK$-model
$y^2 = f(x)$ for $E$ there is at least one $\KK$-rational root of
the cubic $f$. 
The discriminant of $f$ is a square in $\KK^\times$, 
so $f$ splits over $\KK$ and hence 
$E[2]$ is $\KK$-split.  

Let $\KK_{\rm{sep}}/\KK$ be a separable closure and let
$\Gamma \subseteq
{\rm{GL}}_2(\ZZ_2)$
be the image of the 
2-adic representation 
$$\rho_{E,2}:{\rm{Gal}}(\KK_{\rm{sep}}/\KK) \rightarrow
{\rm{Aut}}(\invlim E[2^n](\KK_{\rm{sep}})) \simeq
{\rm{GL}}_2(\ZZ_2)$$
attached to $E$.  Since $E[2]$ is $\KK$-split, 
$\Gamma$ has trivial reduction modulo 2. 
Thus, the Galois action must be pro-$2$,
and hence tame.  Since $E$ has potentially good reduction
and $\OO_{\KK}$ is strictly henselian, 
$\Gamma$ must be 
a finite cyclic 2-group.  Pick $\gamma \in \Gamma$, 
so $\gamma = 1 + 2x$ with
$x \in {\rm{M}}_2(\ZZ_2)$.   Since
the 2-adic cyclotomic character over $\KK$ is trivial,
$$1 = \det(\gamma) = 1 + 2 {\rm{Tr}}(x) + 4 \det(x) = 
-1 + {\rm{Tr}}(\gamma) + 4 \det(x).$$
In particular, ${\rm{Tr}}(\gamma) \equiv 2 \bmod 4$. 
Elements in ${\rm{GL}}_2(\QQ_2)$ with order 4 have characteristic
polynomial $X^2 + 1$ and thus
have trace 0.  This shows that $\Gamma$
cannot contain elements of order 4, so $\Gamma$ 
is either trivial or has order 2.  
Hence, if $\KK' \subseteq \KK_{\rm{sep}}$ is the splitting
field for $\rho_{E,2}$ then $[\KK':\KK]$ divides 2.

The quadratic twist $E'$ of $E$ by $\KK'/\KK$
must have trivial 2-adic
representation, and hence it has good reduction.
Since $\OO_{\KK}$ is strictly henselian and
$[\KK':\KK]$ divides 2, 
$$v(\Delta(E')) \equiv v(\Delta(E)) + 12/[\KK':\KK] \equiv
12/[\KK':\KK] \bmod 12$$
and $v(\Delta(E')) \equiv 0 \bmod 12$ (since $E'$ has
good reduction).  Thus, $[\KK':\KK] = 1$, 
so $E \simeq E'$ has good reduction.
\end{proof}

\begin{lemma}\label{bad}
Assume $\chr(k) \ne 2$ and let $E$ be 
an elliptic curve over $\KK$ with potentially multiplicative reduction.
The parameters $c_4$ and $c_6$ are non-zero, and 
there is multiplicative reduction if and only if
$v(c_4) \equiv 0 \bmod 4$. When 
$\OO_{\KK}$ is complete, there is split
multiplicative reduction if and only if $-c_6$ is a square in $\KK^\times$.
\end{lemma}

\begin{proof}
By hypothesis, $j$ is non-integral.  
As is well-known, $c_4$ and $c_6$ are non-zero when
$j \ne 0, 1728$.
The formation of the N\'eron model
commutes with base change to the completion, so we may
suppose $\OO_{\KK}$ is complete.
Since $2 \in \KK^{\times}$, 
the quadratic extensions of $\KK$
are classified by $\KK^{\times}/(\KK^{\times})^2$,
and since $2 \in \OO_{\KK}^{\times}$, 
the unramified quadratic extensions of $\KK$
are classified by the unit classes
(modulo unit squares).  For $a \in \KK^{\times}$,
let $E^{(a)}$ denote the quadratic twist of $E$ by
the non-trivial character of ${\rm{Gal}}(\KK(\sqrt{a})/\KK)$.  Clearly
$$c_4(E^{(a)}) \equiv a^2 c_4(E) \bmod (\KK^{\times})^4,\,\,\,
c_6(E^{(a)}) \equiv a^3 c_6(E) \bmod (\KK^{\times})^6.$$

By the theory of Tate models, there is a unique class
$u = u(E)$ modulo
$(\KK^{\times})^2$ such that $E^{(u)}$ has split multiplicative reduction. 
Moreover, 
\begin{itemize}
\item $E$ has multiplicative reduction 
if and only if $u$ is a unit class in $\KK^\times/(\KK^\times)^2$, 
\item $E$ has split 
multiplicative reduction if and only if $u$ is trivial in 
$\KK^\times/(\KK^\times)^2$. 
\end{itemize}
A direct calculation with
Tate models shows $c_4(E^{(u)}) \in (\KK^{\times})^4$ 
and $-c_6(E^{(u)}) \in (\KK^{\times})^2$ 
(since $2 \in \OO_{\KK}^{\times}$).
We conclude that
$$-c_6(E) \equiv u \bmod (\KK^{\times})^2,$$
so the reduction is split multiplicative if and only if
$-c_6(E)$ is a square in $\KK^{\times}$.
Also, $u \in \KK^{\times}/(\KK^{\times})^2$ is a unit class if and only if
$v(u) \equiv 0 \bmod 2$, or equivalently
$v(u^2) \equiv 0 \bmod 4$. Therefore, since 
$$v(c_4(E)) \equiv v(c_4(E^{(u)})) - v(u^2) \equiv -v(u^2) \bmod 4,$$
the reduction is multiplicative if and only if $v(c_4(E)) \equiv 0 \bmod 4$.
\end{proof}

Since we are interested in working with
elliptic curves that are not in characteristic 2
and have a non-zero rational 2-torsion point, the shape of a
Weierstrass model can be taken to be
\begin{equation}\label{Weqn}
E:y^2 = x^3 + Ax^2 + Bx.
\end{equation}

The discriminant $\Delta$ and parameters $c_4$ and $c_6$ of
such a model are given by the following
formulas \cite[p.~46]{silverman1}:
\begin{equation}\label{gen}
\Delta = 16B^2(A^2-4B), \ \
c_4 = 16(A^2-3B), \ \ c_6 = -32A(2A^2-9B). \ \
\end{equation}
For $P = (x,y) \in E - E[2]$, the point $[2]P$
has coordinates given by \cite[pp.~58--59]{silverman1}: 
\begin{equation}\label{double}
[2]P = 
\left(\left(\frac{x^2 - B}{2y}\right)^2, 
-\frac{3x^2+2Ax+B}{2y}\left(\frac{x^2 - B}{2y}\right)^2 
+ \frac{x^3 - Bx}{2y}\right).
\end{equation}

We set $A = T$ and $B = -T^3$ in (\ref{Weqn}), giving the elliptic curve
\begin{equation}\label{Eeqn}
E_T : y^2 = x^3 + Tx^2 - T^3x
\end{equation}
over $\FF_p(T)$ with $p \ne 2$. 
By (\ref{gen}), the discriminant and $j$-invariant of (\ref{Eeqn}) are
\begin{equation}\label{Delta}
\Delta = 16T^8(1+4T), \ \ \ 
j = \frac{c_4^3}{\Delta} = \frac{256(1+3T)^3}{T^2(1+4T)}, 
\end{equation}
and the parameters $c_4$ and $c_6$ are
\begin{equation}\label{c4c6}
c_4 = 16T^2(1+3T), \ \ \ c_6 = -32T^3(2+9T).
\end{equation}
For each $t \in F = \kappa(u)$ with $t \not\in \kappa$, the
Weierstrass model
\begin{equation}\label{Wmain}
E_t : y^2 = x^3 + tx^2-t^3x
\end{equation}
defines an elliptic curve over $F$.
(If $t \in \kappa - \{0,-1/4\}$ then $E_t$ is also 
an elliptic curve over $F$, but assuming $t \not\in \kappa$ 
will avoid some unnecessary complications.)  

\begin{theorem}\label{redtype}
Fix $t \in F = \kappa(u)$ with $t \not\in \kappa$.  Let $v$ be a place on $F$.
The reduction type of $E_t$ at $v$ is as in Table $\ref{redtable}$.

\begin{table}
\begin{center}
\begin{tabular}{|cc|}\hline
 $v(t)$          & {\rm Reduction Type} \\ \hline
$>0$, {\rm even} & {\rm multiplicative} \\
$>0$, {\rm odd} & $(${\rm{pot. mult.}}$)$ {\rm additive} \\
$< 0$, $\equiv 0 \bmod 4$ & {\rm good}\\
$< 0$, $\not\equiv 0 \bmod 4$ & $(${\rm{pot. good}}$)$ {\rm{additive}} \\
$=0$, $v(1+4t) = 0$ & {\rm{good}}\\
$=0$, $v(1+4t) > 0$ & {\rm{multiplicative}}\\
\hline
\end{tabular}
\caption{Reduction types for $E_t$, $t \in F - \kappa$}\label{redtable}
\end{center}
\end{table}
\end{theorem}

\begin{proof}
Specializing (\ref{Delta}) and (\ref{c4c6}),  the parameters of
$E_t$ are
\begin{equation}\label{Dceq}
\Delta = 16t^8(1+4t), \ \ j = \frac{256(1+3t)^3}{t^2(1+4t)}, \ \ 
c_4 = 16t^2(1+3t), \ \ c_6 = -32t^3(2+9t).
\end{equation}
(We write $\Delta$ instead of $\Delta|_{T=t}$, and likewise for the 
other parameters.) 
None of the parameters in
(\ref{Dceq}) is 0, since $t \not\in \kappa$. 

If $v(t) > 0$, then
$$
v(\Delta) = 8v(t), \ \ v(c_4) = 2v(t), 
\ \ v(j) = -2v(t) < 0,
$$
so there is potentially multiplicative reduction. Using
Lemma \ref{bad}, there is 
multiplicative reduction when $v(t)$ is even, and there is 
additive reduction when 
$v(t)$ is odd.

If $v(t) < 0$ and $\chr(\kappa) > 3$, then 
$$
v(\Delta) = 9v(t), \ \ v(c_4) = 3v(t), \ \ v(j) = 0,
$$
so there is potentially good reduction.  
If $v(t) < 0$ and $\chr(\kappa) = 3$, then
$$
v(\Delta) = 9v(t), \ \ v(c_4) = 2v(t), \ \ v(j) = -3v(t) > 0,
$$
so again there is potentially good reduction.
Using Lemma \ref{good} in 
both cases, there is good reduction when $v(t) \equiv 0 \bmod 4$ and
there is additive reduction otherwise. 

Finally, suppose $v(t) = 0$, so
$$
v(\Delta) = v(1+4t), \ \ \
v(c_4) =  v(1+3t).
$$
Both $1+4t$ and $1+3t$ have non-negative valuation at $v$, and the valuations
are not both positive.  
If $v(1+4t) = 0$ then $v(j) = 3v(c_4) \geq 0$, so there is good
reduction (by Lemma \ref{good}).
If $v(1+4t) > 0$ then $v(c_4) = 0$, so $v(j) = -v(\Delta) < 0$.
This implies (by Lemma \ref{bad}) that there is 
multiplicative reduction at $v$ in such cases.
\end{proof}

Now we turn to the Mordell--Weil group of the 
generic fiber, $E_T(\FF_p(T))$.  As before, $p \not= 2$.
Two obvious non-zero rational points are $(0,0)$ and
$Q = (-T,T^2)$. 
(There is another obvious non-zero rational point, $(T^2,T^3)$, but 
this is $(0,0)+Q$.)  We will prove that $Q$ has infinite 
order, so ${\rm rank}(E_T(\FF_p(T))) \geq 1$. 

For elliptic curves over $\QQ$, explicit rational points 
are usually checked to be non-torsion by 
the Nagell--Lutz integrality criterion. 
This criterion is really a collection of local criteria
over $\ZZ_{(p)}$ for all primes $p$.
We need an analogue for discrete valuation rings
with positive characteristic. 
Here is a version over arbitrary discrete valuation rings.

\begin{theorem}\label{nlutzFT}
Let $R$ be a discrete valuation
ring with residue field
$k$ of characteristic $p \ge 0$, and let $\KK$ be its
fraction field. Let $E_{/\KK}$ be an elliptic curve,
and let $P \in E(\KK)$ be a non-zero torsion point.

If there exists a Weierstrass model of $E$ over $R$ such that 
one of the affine coordinates of $P$ does not lie in $R$, then 
the scheme-theoretic closure
of $\langle P \rangle \subseteq E(\KK)$ in the N\'eron model
of $E$ over $R$ is a finite flat
local $R$-group.  In particular, $p > 0$ and $P$ has $p$-power order.
If in addition $\chr(\KK) = p$, then 
$E_{/\KK}$ has potentially supersingular reduction and
$j(E) \in \KK$ is a $p$th power.
\end{theorem}

It follows from the Oort--Tate classification
and Cartier duality that the only example of a non-trivial finite flat
local group scheme over $\ZZ_{(p)}^{\rm{sh}}$ with $p$-power order 
and cyclic constant generic fiber is $\mu_2$ for $p = 2$.
Thus,  Theorem \ref{nlutzFT}
recovers the integrality of non-trivial torsion points
on Weierstrass $\ZZ$-models of
the form $y^2 = f(x)$ (for which non-zero $2$-torsion
points must have the form $(x_0,0)$ with $x_0 \in \ZZ$).

\begin{proof}
Let $W \subseteq \PP^2_R$ be the chosen
Weierstrass $R$-model for $E$ (so
there is a chosen isomorphism $W_{\KK} \simeq E$
as pointed curves over $\KK$).  Let
$W^{\rm{sm}} \subseteq W$ be the open $R$-smooth locus of $W$, and let
$\varepsilon \in W(R) = W(\KK) = E(\KK)$ be the 
section $[0,1,0]$, so $\varepsilon \in W^{\rm{sm}}(R)$.
Since $P$ viewed 
as a point 
$$\widetilde{P} \in W(\KK) - \{[0,1,0]\} \subseteq
{\mathbf{A}}^2(\KK) = \KK \times \KK$$
is assumed to have at least one coordinate not in $R$,
as a point of $\PP^2(\KK) = \PP^2(R)$ its reduction in $\PP^2(k)$ cannot
lie in ${\mathbf{A}}^2_k$.  Therefore, the reduction 
must lie on the line
at infinity.  However, by the theory of Weierstrass
models we know that $W_k$ has $\varepsilon_k$
as its unique point on this line, so the reduction
of $\widetilde{P}$ is $\varepsilon_k$.
Since $\varepsilon_k \in W_k^{\rm{sm}}$, we conclude
that $\widetilde{P} \in
W^{\rm{sm}}(R)$.  

Let $\mathscr{E}$ be the N\'eron model of $E$ over $R$. 
By the N\'eron mapping property, there is a unique
map $W^{\rm{sm}} \rightarrow \mathscr{E}$ over $R$
extending the identification of $\KK$-fibers with $E$.
This map carries $\varepsilon$ to the identity element
in $\mathscr{E}(R)$, so the image of $\widetilde{P}$ in
$\mathscr{E}(R)$ reduces to the
identity in $\mathscr{E}_k$.  In other words, under
the equality $E(\KK) = \mathscr{E}(R)$, the reduction of
$P$ in the closed fiber $\mathscr{E}_k$ of the N\'eron model
must be the identity. 

Let $N > 1$ be the order of $P$.  By the N\'eron mapping property,
$P$ defines a map of $R$-groups
$\ZZ/N\ZZ \rightarrow \mathscr{E}$
that is a closed immersion on the generic fiber.
Since the $R$-group $\ZZ/N\ZZ$ is proper and
the target $\mathscr{E}$ is separated 
over the Dedekind domain $R$, 
the scheme-theoretic image of this map is
a finite flat $R$-subgroup $G \hookrightarrow \mathscr{E}$
with order $N$ and constant generic fiber; this
must be the closure of $\langle P \rangle$.
The closed fiber of $G$ must be infinitesimal
since $P$ has reduction equal to the identity.
This forces $G$ to be local, so 
the characteristic $p$ of $k$ must be positive
and the order $N$ of $G$ must be a power of $p$.

Now assume ${\rm{char}}(\KK) = p$.
We must prove that $j(E)$ is a $p$th power in $\KK$ and that
$E$ has potentially supersingular reduction.  
To prove that $j(E)$ is a $p$th power in $\KK$
when $E(\KK)$ contains a non-trivial point with $p$-power
order, we use
the classical fact that
if $L$ is a field with characteristic $p > 0$
and $E$ is an elliptic curve over $L$ such that there exists
an \'etale subgroup $\Gamma \subseteq E$ 
with order $p^n$ for some $n \ge 1$ (that is, $E$ is ordinary and 
the connected-\'etale sequence of $E[p^n]$ splits over $L$), 
then $j(E) \in L$ is a $p^n$th power in $L$. To prove this fact, 
let $E' = E/\Gamma$, so the isogeny
$E' \rightarrow E$ that is dual to the
projection $E \rightarrow E'$ has
kernel that is Cartier-dual to $\Gamma$.
This kernel is therefore multiplicative with $p$-power
order, so it is 
infinitesimal.  Since $E'$ is a 1-dimensional
abelian variety in characteristic $p$, it
contains a unique infinitesimal subgroup of order $p^n$.
The relative $n$-fold Frobenius $E' \rightarrow
{E'}^{(p^n)}$ has this subgroup as its kernel, so 
the two quotients $E$ and ${E'}^{(p^n)}$ of $E'$
are $L$-isomorphic as quotients of $E'$. In particular
$j(E) = j({E'}^{(p^n)}) = j(E')^{p^n}$ in $L$.

Finally, returning to our initial situation, we must
show that $E$ has potentially supersingular reduction if
$\KK$ has characteristic $p > 0$.  The assumptions on $R$ and 
on the coordinates of $P$ are unaffected
by replacing $\KK$ with a finite separable extension $\KK'$ and
replacing $R$ with a maximal-adic localization of its
integral closure in $\KK'$,
so we may assume that $E$ has semistable reduction.  It must be
proved that the N\'eron model $\mathscr{E}$ in this 
case has fibral identity component
$\mathscr{E}_k^0$ that is a supersingular elliptic curve.
Assume to the contrary, so $\mathscr{E}_k^0$ is either
a torus or an ordinary elliptic curve; we seek a contradiction.
In either case, the finite local subgroups of
$\mathscr{E}_k^0$ are multiplicative.  Hence,
if we construct $G$ as we did above (the scheme-theoretic closure 
of $\langle P \rangle$ in $\mathscr{E}$) then the infinitesimal closed
fiber $G_k \hookrightarrow \mathscr{E}_k$ must lie
in $\mathscr{E}_k^0$, so $G_k$ 
is multiplicative with $p$-power order.  Since the generic
fiber $G_{\KK}$ is constant, we conclude
that the Cartier dual $G^{\vee}$ has multiplicative
generic fiber.  However, 
$G^{\vee}$ is finite and flat over $R$ with special fiber
$G_k^{\vee}$ that is \'etale, so $G^{\vee}$ is $R$-\'etale.
This forces $G^{\vee}_{\KK}$ to be both multiplicative and \'etale,
but a non-zero \'etale 
$\KK$-group with $p$-power order cannot be multiplicative when 
$\KK$ has characteristic $p$, so we have reached a contradiction.
\end{proof}

\begin{corollary}\label{djp}
With notation as in Theorem $\ref{nlutzFT}$,
if ${\rm{char}}(\KK) = p > 0$ and $j(E) \in \KK$ is not
a $p$th power, then
non-zero torsion points in $E(\KK)$ have integral
coordinates with respect to any Weierstrass $R$-model of $E$.
\end{corollary}

As an application of Corollary \ref{djp}, we have:

\begin{corollary}\label{rankup}
Assume $p \ne 2$.  The $\FF_p(T)$-rational point
$Q = (-T,T^2)$
on the elliptic curve $E_T$ in $(\ref{Eeqn})$ 
has infinite order.  

In particular, for any field $L$ 
of characteristic $p$ and any $t \in L$
that is transcendental 
over $\FF_p$, 
the specialization $Q_t \in E_t(L)$ that is obtained by 
sending $\FF_p(T)$ into $L$ by $T \mapsto t$
is a point of infinite order.
\end{corollary}

\begin{proof} 
The second claim follows from the first since the 
field extension
$\FF_p(T) 
\rightarrow L$ defined by $T \mapsto t$ induces an injection of groups 
$E_T(\FF_p(T)) \rightarrow E_t(L)$.  

To see that $Q$ has infinite order in $E_T(\FF_p(T))$, first note the 
$j$-invariant of $E_T$, as given in (\ref{Delta}), 
is not a $p$th power in $\FF_p(T)$.
Therefore, by Corollary \ref{djp}, an $\FF_p(T)$-rational point on 
$E_T$ has infinite order 
provided that, using the Weierstrass model (\ref{Eeqn}) for $E_T$, 
some non-zero 
multiple of the point has an $x$- or $y$-coordinate that is non-integral at a 
finite place on $\FF_p(T)$.  (The Weierstrass model 
(\ref{Eeqn}) is integral away from $\infty$.) 

Since $x(Q)$ and $y(Q)$ are 
integral away from $\infty$ and $y(Q) \ne 0$, we double $Q$.
By (\ref{double}), 
$$
[2](Q) = 
\left(
\left(\frac{T+1}{2}\right)^2,\left(\frac{(T+1)(T^2-4T-1)}{8}\right)
\right).
$$
Thus, $x([2]Q)$ and $y([2]Q)$ are integral away from $\infty$ 
and $y([2]Q) \ne 0$, so we double again and find
$$
x([4](Q)) = \left(\frac{(T+1)^4 + 16T^3}{4(T+1)(T^2-4T-1)}\right)^2.
$$
This is non-integral at the place $T + 1$, so we are done.
\end{proof}

\begin{remark} 
By the Lang--N\'eron theorem, the group 
$E_T(\overline{\FF}_p(T))$ is finitely generated.
This group has rank at least 1, since we have exhibited an explicit
element with infinite order. 
Theorem \ref{mainthm} implies that $E_T(L)$ has rank 1
for certain extensions $L$ of
$\overline{\FF}_p(T)$ with transcendence degree 2 over 
$\FF_p$,
so {\em a posteriori} we conclude that $E_T(\overline{\FF}_p(T))$ has rank 1.
Presumably the proof
of Theorem \ref{mainthm} can be modified 
to give a direct proof that $E_T(\overline{\FF}_p(T))$ has rank 1, 
without requiring the use of such auxiliary fields $L$.
\end{remark}

\section{Root numbers}\label{rtlocnum}

For $E_T$ as in (\ref{Eeqn}), 
we will compute the local root numbers $W_v(E_t)$
for $t \in \kappa(u)$ with $t \not\in \kappa$. 
Let us first collect a general list of local root number formulas.
This is well-known for residue characteristic $p \ne 2, 3$,
but we include some cases with $p = 3$.

\begin{theorem}\label{localrtnum}
Let $\KK$ be a local field, with finite residue field of characteristic
$p \not= 2$ and normalized valuation
$v:\KK^{\times} \rightarrow \ZZ$.  
Let $\chi_\KK$ be the quadratic character of
the residue field of $\KK$, and 
let $E$ be an elliptic curve over $\KK$.

\begin{enumerate}
\item Assume $E$ has potentially good reduction,
and if $p = 3$ then assume
$E(\KK)[2] \ne O$. Define $e = 12/\gcd(v(\Delta),12)$.
We have $e \in \{1,2,3,4,6\}$, 
with $3\nmid e$ when $p = 3$, and
the local root number $W_\KK(E)$
can be computed by the following formulas:
$$
W_\KK(E) =
\begin{cases}
1 & \text{ if } e = 1, \\
\chi_\KK(-1) & \text{ if } e = 2 \text{ or } 6,\\
\chi_\KK(-3) & \text{ if } e = 3,\\
\chi_\KK(-2) & \text{ if } e = 4.
\end{cases}
$$

\item Suppose $E$ has potentially multiplicative reduction.  If the
reduction is additive then $W_\KK(E) = \chi_\KK(-1)$.
If the reduction is multiplicative and $c_6 = c_6(E)$ is computed 
using any Weierstrass $\KK$-model of $E$, then
$W_\KK(E) = -1$ when $-c_6$ is a square in $\KK^\times$ and
$W_\KK(E) = 1$ when $-c_6$ is not a square in $\KK^\times$.
\end{enumerate}
\end{theorem}

\begin{proof}
We first address the properties of $e$ in the cases
with potentially good reduction.
By the method of proof of Lemma \ref{good}(2), if $E$ has potentially
good reduction then it acquires good reduction over
a quadratic extension of a splitting field for $E[2]$.
This splitting field
is a tame Galois extension with degree dividing 6, so
Lemma \ref{good} implies that $e$ must divide $12$
and moreover that if $p = 3$ then $3\nmid e$. 
The tameness and Lemma \ref{good}
ensure that $e$ is the order of the image of
inertia in the $\ell$-adic representation for $E$ (any $\ell \ne p$).
The cyclicity of tame inertia therefore rules out
the possibility $e = 12$, since there are
infinitely many rational primes $\ell > 3$
for which the 12th cyclotomic polynomial
$\Phi_{12}$ has no quadratic factors over
$\QQ_{\ell}$.

Before we treat the general case, let us consider
the special case $\KK = \QQ_p$ with $p \not= 2$.  In this case, 
the proposed formulas are proved by Rohrlich for $p > 3$
in \cite[Prop. 2]{rohrlich1}
when the reduction is potentially good
and (using Lemma \ref{bad}) in \cite[Prop.~3]{rohrlich1} when
the reduction is potentially multiplicative.
(Also see \cite[Prop.~3]{rohrlich2} for further
discussion in the multiplicative case.)  
By Lemma \ref{good}(2) and Lemma \ref{bad},
the proofs of
\cite[Prop.~2,~3]{rohrlich1} work in our cases when $p = 3$.
(The additional 2-torsion hypothesis in potentially good
reduction cases for $p = 3$ avoids wild ramification.) 

Rohrlich's proofs in \cite{rohrlich1} and
\cite{rohrlich2} are representation-theoretic and
rest on papers of Deligne~\cite{deligne}
and Tate~\cite{tate2} that are valid for
local fields with any residual (or generic)
characteristic.  Thus, these proofs 
carry over to the general case (with residue characteristic $\not= 2$,
and with a non-trivial rational 2-torsion point in
potentially good cases when $p = 3$).
The ``$p$'' in
 most of the arguments in \cite{rohrlich1} and \cite{rohrlich2}
should be replaced with the size of the residue field of $\KK$, say $q$,
and the Legendre symbol $(\frac{a}{p})$ should be replaced
with the Kronecker symbol
$(\frac{a}{q})$. (Note that when $q$ is an odd prime
prime power and $a \in \ZZ$ is prime to $q$, $(\frac{a}{q}) =
\chi_\KK(a)$.)  A general discussion
in the context of local and global fields of characteristic 0
may also be found in \cite{rohrlich3}.
\end{proof}

Using Theorems \ref{redtype} and \ref{localrtnum}, 
we now compute the local root numbers $W_v(E_t)$ for every $t \in
F = \kappa(u)$ with $t \not\in \kappa$.
Let $\chi$ be the quadratic character of $\kappa$ and
let $\chi_v$ be the quadratic character of the residue field at $v$.
For $a \in \kap$, we have $\chi_v(a) = \chi(a)^{\deg v}$, 
where $\deg v$ is the degree of the residue field of $v$ over 
$\kappa$. 
(Thus, $\chi_v = \chi$ when $\deg v = 1$.) 
Table \ref{roottable} summarizes the results, and we will see 
why the first row is undesirable.

\begin{table}
\begin{center}
\begin{tabular}{|cc|}\hline
 $v(t)$               & $W_v(E_t)$ \\ \hline
$> 0$, even                & hard to use\\ 
$> 0$, odd            & $\chi_v(-1)$ \\
$< 0$, even            & $\chi_v(-1)^{v(t)/2}$ \\
$< 0$, odd            & $\chi_v(-2)$ \\
$=0$, $v(1+4t) = 0$   & 1 \\
$=0$, $v(1+4t) > 0$  & $-\chi_v(2)$ \\
\hline
\end{tabular}
\caption{Local root numbers on $E_t$, $t \in 
F - \kappa$}\label{roottable}
\end{center}
\end{table}

The second, third, fourth, and fifth rows are cases of additive or 
good reduction (by Table \ref{redtable} in Theorem \ref{redtype}), 
and these are left to the reader to check via 
Theorem \ref{localrtnum}.
(The third row is the union of two cases from Table \ref{redtable}
with $v(t) < 0$, namely $v(t) \equiv 0, 2 \bmod 4$.  These two cases 
are checked separately.)
It remains to compute $W_v(E_t)$ in two cases:
(i) $v(t)$ is positive and even, and
(ii) $v(1+4t) > 0$.  Both are cases of multiplicative reduction, so 
Theorem \ref{localrtnum} requires us to check if $-c_6(E_t)$ 
is a square in the multiplicative group $F_v^{\times}$ of the completion 
of $F$ at $v$.  Let us first introduce some convenient notation.

\begin{definition}\label{funnydef}
Let $L$ be a field.
For $x, y \in L$, write $x \sim y$ when $x = y z^2$
for some $z \in L^{\times}$.
\end{definition}

When $v(t)$ is positive and even, in $F_v^\times$ we compute from 
(\ref{c4c6}) at $T = t$ that 
\begin{eqnarray*}
-c_6 & = &
32t^3(2+9t) \\
& \sim & 2t(2+9t) \\
& \sim & t \ \ \  \text{ since } v(t) > 0. \\
\end{eqnarray*}
Thus, by Theorem \ref{localrtnum}(2), $W_v(E_t) = -1$ if 
$t$ is a square in $F_v^\times$ and $W_v(E_t) = 1$ otherwise. 

The last case is $v(1+4t) > 0$. In $F_v^\times$, 
\begin{eqnarray*}
-c_6 & = & 32t^3(2+9t) \\
& \sim & 2t(2+9t) \\
& \sim & 2t^2 \ \ \  \text{ since } v(1+4t) > 0 \\
& \sim & 2.
\end{eqnarray*}
Thus, by Theorem \ref{localrtnum}(2), the last entry 
in Table \ref{roottable} is confirmed:
\begin{equation}\label{WvEt}
v(1+4t) > 0 \Longrightarrow W_v(E_t) = -\chi_v(2).
\end{equation}

The $E_t$'s do not have easily manageable global root numbers.
There are two main problems.
First, the local root number in the first row 
of Table \ref{roottable} depends on whether or not $t$ is a 
square in $F_v^\times$, and that is not something we can easily 
control.  Second, the last row in
Table \ref{roottable} introduces {\em systematic}
minus signs.
To appreciate the nature of these difficulties, and how 
we can avoid them by a change of variables that is peculiar to 
characteristic $p$, let us first try to eliminate
the difficulties in the first row of
Table \ref{roottable} by forcing ``$t$''
to be a square:
we study the elliptic curve $E_{T^2}$. 
Table \ref{roottable} is easily translated into this context, 
and the results are collected in Table \ref{roottable2}. 
The systematic minus signs in the first and last
rows of Table \ref{roottable2} will cause serious problems.

\begin{table}
\begin{center}
\begin{tabular}{|cc|}\hline
 $v(t)$              & $W_v(E_{t^2})$ \\ \hline
$> 0$             & $-1$\\ 
$< 0$             & $\chi_v(-1)^{v(t)}$ \\
$=0$, $v(1+4t^2) = 0$   & 1 \\
$=0$, $v(1+4t^2) > 0$   & $-\chi_v(2)$ \\
\hline
\end{tabular}
\caption{Local root numbers on $E_{t^2}$, $t \in 
F - \kappa$}\label{roottable2}
\end{center}
\end{table}

Write $t = g_1/g_2$, where $g_1, g_2 \in \kappa[u]$ 
are non-zero and relatively prime.  The
 product of the local 
root numbers $W_v(E_{t^2})$ over all $v$ 
yields the global root number formula 
\begin{eqnarray*}
W(E_{t^2}) 
& = & W_{\infty}(E_{t^2})\prod_{\Atop{v \not= \infty}{v(t) > 0}} 
(-1) \cdot \prod_{\Atop{v \not= \infty}{v(t) < 0}}\chi_v(-1)^{v(t)} 
\cdot \prod_{\Atop{v \not= \infty}{v(1+4t^2) > 0}}(-\chi_v(2)) \\
& = & W_{\infty}(E_{t^2})\cdot (-1)^{\#\{\pi : \pi|g_1\}}
\cdot \chi(-1)^{\sum_{\pi|g_2}(\deg \pi)\ord_\pi(g_2)}
\cdot \prod_{\pi|(4g_1^2+g_2^2)} (-\chi(2)^{\deg \pi})  \\
& = & W_{\infty}(E_{t^2})\cdot (-1)^{\#\{\pi : \pi|g_1\} +
\#\{\pi : \pi|(4g_1^2+g_2^2)\}}
\cdot \chi(-1)^{\deg g_2}
\cdot \chi(2)^{\sum_{\pi|(4g_1^2+g_2^2)} \deg \pi}, 
\end{eqnarray*}
where $\pi$ runs over monic irreducibles in $\kappa[u]$. 
This formula is unwieldy because we
cannot control the parity of the number of irreducible factors of $g_1$
or $4g_1^2 + g_2^2$ 
as $t$ varies.  
We also cannot control the parity
of $\sum_{\pi|(4g_1^2 + g_2^2)} \deg \pi$
because it is hard to determine when
$4g_1^2 + g_2^2$ is separable, though 
this second problem could be eliminated if we only consider
$\kappa$ in which $\chi(2) = 1$.
Studying $E_{t^2}$ is not helping us to get 
constant global root numbers at most $t$ (as is essentially necessary
in any example of elevated rank).

Instead of merely simplifying the first row of Table \ref{roottable} by 
replacing $t$ with $t^2$ in $E_t$, we need to 
eliminate the intervention of the 
first row of Table \ref{roottable}.
To accomplish this, we will introduce a change 
of variables in $t$ such that the numerator is  
 always squarefree, and thus in particular 
never has positive even valuation at places of $F$. 
We also need to acquire
control over the product of minus signs contributed from 
the last row of Table \ref{roottable}, and this will be
achieved by arguments that are peculiar to positive characteristic. 

A ``squarefree'' change of variables is impossible 
in characteristic 0, but the $p$th power map provides a mechanism 
to find such a change of variables
in characteristic $p$.  The basic idea is that, for all 
$t \in F = \kappa(u)$, $t^p + u$ has a squarefree numerator
and has a pole at $\infty$, and thus, for all places $v$ of $F$, 
$v(t^p+u)$ is never both positive and even.  With this noted, define
\begin{equation}\label{hdef}
h(T) = cT^{2p} + du,
\end{equation}
where $c, d \in \kappa^\times$.
The use of the exponent $2p$ instead of $p$ will create 
a counterexample to Chowla's two-variable conjecture
over $\kappa[u]$ (see the appendices for a discussion of this
conjecture and its relevance to the study of elevated rank); 
in concrete terms, this even exponent will force certain
otherwise unknown non-zero quantities we meet later to be squares.
The role of $c$ and $d$ in $h(T)$ is to provide us 
with the family of examples in Theorem \ref{mainthm} for each
$p \ne 2$, 
rather than just one example for each $p \ne 2$.
(The reader may take $c = d = 1$ throughout.)

For $h(T)$ as in (\ref{hdef}), consider the elliptic curve $E_{h(T)}$ over
$F(T)$, obtained by replacing $T$ with $h(T)$ in (\ref{Eeqn}). 
We can run through all of our previous 
work with $h(t)$ in place of $t$ (rather than $t^2$ in place of $t$), 
and now $t$ can be any element of $F$ since $h(t) \not\in \kappa$ for all $t 
\in F$.  (Table \ref{roottable2} only lists
root numbers in fibers over $t \in F - \kappa$.)
We will compute $W(E_{h(t)})$ for all $t \in F$. 

Write $t = g_1/g_2$, where  
$g_1, g_2 \in \kappa[u]$ are relatively prime with $g_2 \ne 0$, so
\begin{equation}\label{heqn}
h(t) = \frac{cg_1^{2p} + dug_2^{2p}}{g_2^{2p}}.
\end{equation}
Call the numerator and denominator, respectively, $f_1$ and $f_2$:
\begin{equation}\label{f1f2def}
f_1 = cg_1^{2p} + dug_2^{2p}, \ \ \ f_2 = g_2^{2p}.
\end{equation}
Obviously $f_1, f_2 \ne 0$. 
Since $(g_1,g_2) = 1$, clearly $(f_1,f_2) = 1$. 
Moreover, since $f_1' = dg_2^{2p} = df_2$, $f_1$ is squarefree.
Thus, for all finite places $v$ of $F$,
$v(h(t))$ is never both positive and even at $v$. 
Since 
\begin{equation}\label{ordinfh}
\ord_\infty(h(t)) = \ord_{\infty}(ct^{2p} + du) = 
\begin{cases}
-1 & \text{ if } \ord_{\infty}(t) \ge 0, \\
2p\ord_{\infty}(t) & \text{ if } \ord_{\infty}(t) < 0, 
\end{cases}
\end{equation}
we see $h(t)$ has a pole at $\infty$ for 
every $t \in F$. 

We begin computing local root numbers
for $E_{h(t)}$ over $F$ by starting
with the place at $\infty$, 
where $\chi_\infty = \chi$. 
Using (\ref{ordinfh}) and 
Table \ref{roottable} (with $h(t)$ in place of $t$), 
\begin{equation}\label{Winf}
W_\infty(E_{h(t)}) =
\begin{cases}
\chi(-2), & \text{ if } \ord_{\infty}(t) \ge 0,\\
\chi(-1)^{\ord_{\infty}(t)}, & \text{ if } \ord_{\infty}(t) < 0.
\end{cases}
\end{equation}

Now let $v$ be a finite place on $F$.
Since $h(t)$ has a squarefree
numerator, we get from Table \ref{roottable} that
$$
v \not= \infty, \,\,
v(h(t)) > 0 \Longrightarrow v(h(t)) = 1 \Longrightarrow W_v(E_{h(t)}) =
\chi_v(-1) = \chi(-1)^{\deg v}.
$$
Since $h(t)$ 
has a perfect square $g_2^{2p}$
as its denominator, Table \ref{roottable} implies 
$$
v \not= \infty, \,\,
v(h(t)) < 0 \Longrightarrow W_v(E_{h(t)}) = \chi_v(-1)^{v(h(t))/2} = 
\chi(-1)^{(\deg v) \cdot v(g_2)}. 
$$
If $v(h(t)) = 0$, then 
Table \ref{roottable} (with $h(t)$ in place of $t$) tells us 
that if $v(1+4h(t)) = 0$ then $W_v(E_{h(t)}) = 1$, whereas
$$
v(1+4h(t)) > 0 
\Longrightarrow W_v(E_{h(t)}) = -\chi_v(2) = -\chi(2)^{\deg v}.
$$

Combining all of this local information, 
for $t \in F$ 
the global root number $W(E_{h(t)})$ is
\begin{equation}\label{Wform}
W_{\infty}(E_{h(t)})\prod_{\Atop{v \not= \infty}{v(h(t))>0}}\chi(-1)^{\deg v}
\prod_{\Atop{v \not= \infty}{v(h(t)) < 0}}\chi(-1)^{(\deg v)\cdot v(g_2)}
\prod_{\Atop{v \not= \infty}{v(1+4h(t))>0}} (-\chi(2)^{\deg v}).
\end{equation} 
Referring back to Table \ref{redtable} with $h(t)$ in place of $t$, 
the local root numbers for $E_{h(t)}$ at places of multiplicative reduction 
appear in (\ref{Wform}) as the terms in the last product.

Writing
 (\ref{Wform}) in terms of the numerator and denominator of $h(t)$, 
\begin{eqnarray*}
W(E_{h(t)})  & = & 
W_{\infty}(E_{h(t)})\prod_{\pi|f_1}\chi(-1)^{\deg \pi}
\prod_{\pi|f_2}\chi(-1)^{(\deg \pi)\ord_\pi(g_2)}
\prod_{\pi|(4f_1+f_2)} (-\chi(2)^{\deg \pi}) \\
& = & W_{\infty}(E_{h(t)})\prod_{\pi|f_1}\chi(-1)^{\deg \pi}
\prod_{\pi|g_2}\chi(-1)^{(\deg \pi)\ord_\pi(g_2)}
\prod_{\pi|(4f_1+f_2)} (-\chi(2)^{\deg \pi}) \\
& = & W_{\infty}(E_{h(t)})\prod_{\pi|f_1}\chi(-1)^{\deg \pi} 
\cdot \chi(-1)^{\deg g_2} \cdot \prod_{\pi|(4f_1+f_2)} (-\chi(2)^{\deg \pi}).
\end{eqnarray*}
Set 
\begin{equation}\label{fff}
f = 4f_1 + f_2 = 4cg_1^{2p} + (4du+1)g_2^{2p}.
\end{equation}
Since $(f,f') = 1$, $f$ is squarefree.  
We already saw that $f_1$ is squarefree as well, so 
our global root number formula simplifies to 
\begin{equation}\label{Wht}
W(E_{h(t)}) = 
W_\infty(E_{h(t)})\chi(-1)^{\deg f_1}
\chi(-1)^{\deg g_2}\mu(f)\chi(2)^{\deg f},
\end{equation}
where $\mu$ is the M\"obius function on 
$\kappa[u]$ (defined much like its classical counterpart
over $\ZZ$).  

\begin{remark}\label{remlambda}
Let us clarify how this calculation
is analogous to what is seen in work over $\QQ(T)$.
Let the {\em Liouville function} $\lambda$ on $\kappa[u]$ be 
the totally multiplicative function whose value on irreducible
elements is $-1$, so if $f$ is separable
(\ie, is squarefree) in $\kappa[u]$ then
$\mu(f) = \lambda(f)$.  In (\ref{Wht}), 
we therefore have an appearance of $\lambda(f)$ 
as a contribution from local root numbers at 
places of multiplicative reduction.
As is explained in Appendix \ref{avgrtnum}, 
the Liouville function on $\ZZ$ arises in a similar
manner in the study of average root numbers
for elliptic curves over $\QQ(T)$
that have a point of multiplicative reduction on $\PP^1_{\QQ}$. 
Another similarity with the situation in characteristic 0
is that $\lambda$ is being computed on
an element $f \in \kappa[u]$ that is the value 
at $(g_1, g_2)$ of a homogeneous 2-variable polynomial over $\kappa[u]$,
where $g_1$ and $g_2$ are relatively prime. (Consider 
$g_1$ and $g_2$ in (\ref{fff}) as specializations of
independent indeterminates over $\kappa[u]$.)  Compare this
with the appearance of $\lambda(f_{\mathscr{E}}(m,n))$ 
in the discussion at the end of Appendix
\ref{avgrtnum}.
\end{remark}

To 
simplify (\ref{Wht}) further, we compute the 
degrees in the exponents.  This depends on 
the relative sizes of $\deg g_1$ and $\deg g_2$. Let
$$n_1 = \deg g_1, \ \ \ \ n_2 = \deg g_2,$$
with the standard convention $n_1 = -\infty$ when $g_1 = 0$, so 
\begin{equation}\label{f1data}
\deg f_1 = \begin{cases}
2pn_2+1 & \text{ if } n_1 \leq n_2, \\
2pn_1 & \text{ if } n_1 > n_2,
\end{cases}
\ \ \
\deg f_2 =
2pn_2. 
\end{equation}
By inspection, $\deg f_1 > \deg f_2$, 
so 
\begin{equation}\label{fdata}
\deg f = \begin{cases}
2pn_2+1 & \text{ if } n_1 \leq n_2, \\
2pn_1 & \text{ if } n_1 > n_2.
\end{cases}
\end{equation}

Using 
(\ref{Winf}), (\ref{Wht}), 
(\ref{f1data}), and (\ref{fdata}), 
\begin{equation*}
W(E_{h(t)})
 = 
\begin{cases}
\chi(-2)\chi(-1)\chi(-1)^{n_2}\mu(f)
\chi(2) & \text{ if } n_1 \leq n_2, \\
\chi(-1)^{n_2-n_1}\cdot 1 
\cdot \chi(-1)^{n_2} \cdot \mu(f) \cdot 1 & \text{ if } n_1 > n_2,
\end{cases}
\end{equation*}
so
\begin{equation}\label{endglobal}
W(E_{h(t)}) =
\begin{cases}
\chi(-1)^{n_2}\mu(f) & \text{ if } n_1 \leq n_2, \\
\chi(-1)^{n_1}\mu(f) & \text{ if } n_1 > n_2, 
\end{cases}
\end{equation}
where $t = g_1/g_2 \in \kappa(u)$ is
expressed in reduced form and 
$f$ is defined in (\ref{fff}).
(The two cases in (\ref{endglobal}) 
are classified by
the sign of $\ord_\infty(t) = n_2 - n_1$.)

To complete the calculation of $W(E_{h(t)})$ for
$t \in F$, we need to compute $\mu(f)$. 
For this, we use a remarkable fact: 
the M\"obius function in characteristic $p$ is a
more accessible object than its classical counterpart over $\ZZ$. Indeed, 
there is a formula for the M\"obius function on
$\kappa[u]$ {\it other} than its definition.
In particular, the explicit calculation of 
\begin{equation}\label{mufg}
\mu(f) = \mu(4cg_1^{2p} + (4du+1)g_2^{2p}),  
\end{equation}
where $g_1$ and $g_2$ appear through their $p$th powers, 
can be done without factoring.
(Nothing of the sort can be said for classical 
variants such as $\mu_{\ZZ}(m^2+5n^2)$.) 

The alternative M\"obius formula (in Lemma \ref{stick} below)
uses discriminants, so
to avoid any possible confusion on signs and scalar factors,
let us briefly recall how to
define the discriminant of a 
polynomial. For 
any field $K$ and any non-zero polynomial $P = P(u)$ 
in $K[u]$ with degree $n$, 
the {\em discriminant} of $P$ is
\begin{equation}\label{discdef}
\disc_K P :=
(-1)^{n(n-1)/2} \cdot (\lead P)^{n-2} \cdot
\prod_{i=1}^n P'(\gamma_i) \in K,
\end{equation}
where $\gamma_1,\dots,\gamma_n$ are the roots of $P$ (repeated with
multiplicity) in a splitting field and
$\lead P \in K^{\times}$ is the leading coefficient of $P$. 
Obviously ${\rm{disc}}_K(cP) = c^{2n - 2} \cdot 
{\rm{disc}}_K(P)$ for $c \in K^{\times}$.
(In \cite{ccg}, which motivated the work in this section, 
a different definition of the discriminant is used 
that is invariant under $K^{\times}$-scaling of $P$. 
That definition differs from (\ref{discdef})
by an even power of $\lead P$.
Discriminants will only matter up to a non-zero square scaling 
factor for our purposes, 
because of the quadratic character in (\ref{muf}) below, so 
the different discriminants used in \cite{ccg} and here 
are not incompatible for the intended applications.)

\begin{lemma}\label{stick}
Let $\kap$ be a finite field with odd characteristic, and
let $\chi$ be the quadratic character on $\kap$, with $\chi(0) = 0$.
For any non-zero polynomial $P \in \kap[u]$,
\begin{equation}\label{muf}
\mu(P) = (-1)^{\deg P}\chi(\disc_\kap P),
\end{equation}
where $\disc_\kap P$ is the discriminant of the polynomial $P$.
\end{lemma}

\begin{proof}
This 
formula is trivial when $P$ has a multiple factor: both
sides are 0.
When $P$ is separable (that is, squarefree) and has $r$
prime factors, (\ref{muf}) is the same 
as: $\chi(\disc_\kap P) = (-1)^{\deg P - r}$. 
Written this way, (\ref{muf}) appears in \cite[Cor. 1]{swan}.
The properties of finite fields that are
most essential in the proof of \cite[Cor.~1]{swan} are perfectness and
pro-cyclicity of their Galois theory. (The 
only reason to assume
${\rm{char}}(\kappa) \ne 2$ is that
the M\"obius formula can then be given in terms of the quadratic character; 
a formula when
${\rm{char}}(\kappa) = 2$ can 
be found in \cite{ccg} and \cite{swan}, but 
it uses a lift to characteristic 0.  We 
omit this formula since we do not need it.)
\end{proof}

Direct computation of polynomial 
discriminants can often be unwieldy, so
applications of Lemma \ref{stick} usually rest on 
the connection
between discriminants and resultants (see \cite{ccg} 
and 
\cite{swan}); 
such work with resultants requires special care in positive characteristic.
In our specific situation we will be able to 
extract the required
information directly from the definition of the
discriminant, so we will not need
to use resultants.

In (\ref{endglobal}) we have to compute $\mu(f)$ for 
$f = 4cg_1^{2p} + (4du+1)g_2^{2p}$ such that 
$g_2 \ne 0$ and $(g_1, g_2) = 1$.
The peculiar coefficient of $g_2^{2p}$ is an artifact of 
our elliptic curve $E_{h(T)}$.  We shall carry out the 
M\"obius calculation for 
a cleaner expression and then return to $\mu(f)$. 

\begin{lemma}\label{glemma}
Let $\kappa$ be a finite field with characteristic $p \not= 2$. 
Using the convention $\deg(0) = -\infty$, 
for $a, b \in \kappa^\times$ and 
relatively prime $g_1, g_2 \in \kappa[u]$ we have
$$
\mu(ag_1^{2p} + bug_2^{2p}) = 
\begin{cases}
-\chi(-1)^{\deg g_2} & \text{ if } \deg g_1 \leq \deg g_2, \\
\chi(-1)^{\deg g_1} & \text{ if } \deg g_1 > \deg g_2.
\end{cases}
$$
\end{lemma}

\begin{proof}
The cases when $g_1 = 0$ or $g_2 = 0$ are trivial, so we now suppose
both are non-zero. 
Set $g = ag_1^{2p} + bug_2^{2p}$, $n_1 = \deg g_1$, $n_2 = \deg g_2$. 
Since $g' = bg_2^{2p}$ and $(g_1,g_2)= 1$, $g$ is squarefree.
We have (with notation as in Definition \ref{funnydef})
$$
\deg g = 
\begin{cases}
2pn_2+1 & \text{ if } n_1 \leq n_2, \\
2pn_1 & \text{ if } n_1 > n_2,
\end{cases}
\ \ \
\lead g \sim
\begin{cases}
b & \text{ if } n_1 \leq n_2, \\
a & \text{ if } n_1 > n_2.
\end{cases}
$$
Let $n = \deg g$, so in a splitting field the set of
distinct roots of $g$
may be labelled as $\{\gamma_1, \dots, \gamma_n\}$.
Since $g' = b g_2^{2p}$, it
follows that $\prod_i g'(\gamma_i)$ is in $\kappa^{\times}$ and 
may be computed modulo squares:
$$\prod_i g'(\gamma_i) = b^{\deg g} \cdot 
\prod_i g_2(\gamma_i)^{2p} \sim b^{\deg g}$$
because $\prod_i g_2(\gamma_i) \in \kappa^{\times}$. 
Hence, 
$$
{\rm{disc}}_{\kappa}(g) = (-1)^{n(n-1)/2}(\lead g)^{n-2} \cdot
\prod_{i=1}^n g'(\gamma_i) \sim (-1)^{n(n-1)/2}(b \cdot \lead g)^n.$$
Since $b \cdot \lead g \sim b^2$ when $n$ is odd (that is, when
$n_1 \le n_2$), we conclude 
$${\rm{disc}}_{\kappa}(g) \sim (-1)^{n(n-1)/2} \sim
(-1)^{\max(n_1,n_2)}$$
by the formula for
$n$.
Hence, by (\ref{muf}), 
$\mu(g) = (-1)^n \chi({\rm{disc}}_{\kappa}(g)) =
(-1)^n \chi(-1)^{\max(n_1,n_2)}$.
\end{proof}

It is now a simple matter to finish the computation of the global root
number:

\begin{theorem}\label{rootabc}
Let $h(T) = cT^{2p}+du$, where $c, d \in \kappa^\times$.
Let $E_T$ be defined as in $(\ref{Eeqn})$. 
For any $t \in F = \kappa(u)$, the elliptic curve
$E_{h(t)}$ over $F$ satisfies 
\begin{equation}\label{WW}
W(E_{h(t)}) =
\begin{cases}
-1 & \text{ if } \ord_{\infty}(t) \geq 0, \\
1 & \text{ if } \ord_{\infty}(t) < 0.
\end{cases}
\end{equation}
\end{theorem}

\begin{proof}
Write $t = g_1/g_2$ where $g_2 \ne 0$ and $(g_1,g_2) = 1$. 
We may apply Lemma \ref{glemma}
to the polynomial $f = 4cg_1^{2p} + (4du+1)g_2^{2p}$ 
by making the linear change of variables
$u \mapsto u - 1/4d$ that preserves degrees.  This yields
\begin{equation}\label{eq177}
\mu(f) =
\begin{cases}
-\chi(-1)^{n_2}, & \text{ if } n_1 \leq n_2, \\
\chi(-1)^{n_1}, & \text{ if } n_1 > n_2, 
\end{cases}
\end{equation}
where $n_1 = \deg g_1$ and $n_2 = \deg g_2$
(and $n_1 = -\infty$ if $g_1 = 0$). 
Combining (\ref{eq177}) with the global root number formula 
(\ref{endglobal}) yields (\ref{WW}). 
\end{proof}

To force the global root number to be
1, we want only the second case of (\ref{WW}) to occur. 
This can be achieved by a simple trick 
(related 
to (\ref{hdef}), but initially inspired by \cite[p.~57]{kochen}):
replace $t$ with 
$t^2 + u$, which has a pole at $\infty$ for every $t$ in $F = \kappa(u)$.
Thus,
\begin{equation}\label{we1}
W(E_{h(t^2+u)}) = 1
\end{equation}
for {\em every} $t \in F = \PP^1_F(F) - \{\infty\}$.
Since (\ref{mainfamily}) is the Weierstrass model in the definition of 
$E_{h(T^2+u)}$,
we see that (\ref{mainfamily}) is not
as arbitrary as it may have initially
appeared to be.
Combining (\ref{we1}) with Remark \ref{infinity}
settles the root number aspect of Theorem
\ref{mainthm}.

\section{Generic rank bound I. 
Specialization at points of height $0$}\label{rank0}

Write (\ref{mainfamily}) in the form 
\begin{equation}\label{Wut}
\mathscr{E}_{\eta}  : y^2 = x^3 + h(T^2+u)x^2 -(h(T^2+u))^3x, 
\end{equation}
where 
$h(T) = cT^{2p}+du$ 
and $c, d \in \kappa^{\times}$. 
We have shown in \S\ref{rtlocnum} that 
for each $t \in \PP^1(F)$, $\mathscr{E}_t$ is an elliptic curve  
over $F$ with global root number 1. 
The elliptic curve $\mathscr{E}_{\eta}$ over 
$F(T) = \kappa(u,T)$ is obtained from
$E_{T/\FF_p(T)}$ in (\ref{Eeqn}) by replacing $T$ with the element  
$h(T^2+u) \in F(T) = \kappa(u,T)$ that is not in $\kappa$, so
the point $(-T,T^2) \in E_T(\FF_p(T))$ goes
over to the point 
\begin{equation}\label{QdefT}
Q = (-h(T^2+u),(h(T^2+u))^2) \in \mathscr{E}_{\eta}(F(T))
\end{equation}
that has infinite order
(Corollary \ref{rankup}).  For every $t \in F$
the specialization $h(t^2+u) \in \kappa(u)$ is
not in $\kappa$, so
the specialization of $Q$ in 
$\mathscr{E}_t(F)$ must likewise have infinite order
for all $t \in F$. 
Thus, all specializations $\mathscr{E}_t(F)$ 
at $t \in \PP^1(F) - \{\infty\}$ have positive rank.  
This settles the rank aspect of Theorem \ref{mainthm}
for the $F$-rational fibers. (We already noted
in Remark \ref{infinity}
that $\mathscr{E}_{\infty}(F)$ has rank 0.) 

The remainder of this paper is devoted to proving that
the generic Mordell--Weil group 
$\mathscr{E}_{\eta}(F(T))$, which we know has rank at least 1, 
has rank exactly 1.  This
will complete the proof of Theorem \ref{mainthm}.

Since  the cubic polynomial in $x$ given by the Weierstrass
model (\ref{Wut}) defining $\mathscr{E}_{\eta}$ is the product
of $x$ and an irreducible quadratic polynomial in
$F(T)[x]$, the only nontrivial rational 2-torsion is 
the point $(0,0)$. 
Therefore 
\begin{equation}\label{eta2}
\dim_{\FF_2} \mathscr{E}_{\eta}(F(T))/2\cdot \mathscr{E}_{\eta}(F(T)) = 1 +
{\rm{rank}}(\mathscr{E}_{\eta}(F(T))).
\end{equation}
A point of infinite order in $\mathscr{E}_{\eta}(F(T))$
is given in (\ref{QdefT}), so the generic 
rank is 1 if and only if
the dimension in (\ref{eta2}) is at most 2.

Viewing $F(T) = \kappa(u,T)$ as the function field of $\PP^1 \times \PP^1$,
we shall now consider specialization along the $u$-line.
We will specialize at  generic points
of $\{u_0\} \times \PP^1_{\kappa}$
for closed points $u_0 \in \PP^1_{\kappa}$;
these generic points are identified with
the closed points of height 0 on the $u$-line $\PP^1_{\kappa(T)}$
over $\kappa(T)$. 
For such $u_0$, let its residue field be written as 
$
\kappa_0 = \kappa(u_0)$;
this is a finite field and the notation 
$\kappa_0$ will be used constantly in what follows.
If $u_0 \ne \infty$ then we 
also write $u_0$ to denote the image of the indeterminate $u$ under
the quotient map $\kappa[u] \twoheadrightarrow \kappa_0$.

Using (\ref{Delta}) and (\ref{c4c6}), the parameters $\Delta$ and $c_4$ for 
(\ref{Wut}) are given by 
\begin{equation}\label{Dc4wut}
\Delta = 
16(h(T^2+u))^8(1+4h(T^2+u)), \ \ \ 
c_4 = 16(h(T^2+u))^2(1+3h(T^2+u)).
\end{equation}
From the formula for $\Delta$, we see that for all closed points 
$u_0 \in {\mathbf{A}}^1_{\kappa}$, 
the $u_0$-specialization of (\ref{Wut})
is an elliptic curve
over $\kappa_0(T)$.
The elliptic curves $\mathscr{E}_t$ for $t \in \PP^1(F)$ all
live over the fixed global field $F = \kappa(u)$,
but the $u_0$-specializations
$\mathscr{E}_{u_0}$ of $\mathscr{E}_{\eta}$ 
live over the global fields $\kappa_0(T) = \kappa(u_0)(T)$
that vary.  The notation
$\mathscr{E}_{u_0}$ presents no risk of confusion with the notation
${\mathscr E}_{t}$ for specialization at $t \in \PP^1(F)$ because we will 
never again use such $T$-specializations.

Let us briefly describe a natural but ultimately unsuccessful
strategy for using
the $\mathscr{E}_{u_0}$'s to show that the dimension in (\ref{eta2}) 
is at most 2. 
We can prove a ``height 0'' version of Silverman's specialization
theorem for abelian varieties, and from this it follows that
 for all but finitely many height-0 points
$u_0 \in \PP^1_{\kappa(T)}$, 
 the specialization map
\begin{equation}\label{above}
\mathscr{E}_{\eta}(F(T)) \rightarrow \mathscr{E}_{u_0}(\kappa_0(T))
\end{equation}
at $u_0$ is injective.
Thus, it would suffice
to prove ${\rm{rank}}(\mathscr{E}_{u_0}(\kappa_0(T))) \le 1$ for
infinitely many $u_0$.  
For an infinite set of
points $u_0$ (specifically, the ones 
arising from Theorem \ref{phipthm} below), we can
prove ${\rm{rank}}(\mathscr{E}_{u_0}(\kappa_0(T))) \le 3$.
(Switching root number calculations to the $u_0$-side, we also 
can show $W(\mathscr{E}_{u_0}) = -1$. This suggests, but does 
not prove, that $\mathscr{E}_{u_0}(\kappa_0(T))$ has rank 1 or 3.)
For such $u_0$, the
subspace $V_{u_0}$ of everywhere-unramified classes in the
2-Selmer group $S^{[2]}({\mathscr{E}_{u_0}}_{/\kappa_0(T)})$
is 2-dimensional,
and we can show that
${\rm{rank}}(\mathscr{E}_{u_0}(\kappa_0(T))) = 1$ 
(resp. $< 3$) if and only if
the natural map
$V_{u_0} \rightarrow 
\Sha(\mathscr{E}_{u_0})[2]$ is injective (resp. non-zero).
The Cassels--Tate pairing
of the image of a basis of $V_{u_0}$ in $\Sha(\mathscr{E}_{u_0})[2]$
can be calculated by using
a method of Cassels \cite{cassels},
but unfortunately it always turns out to be trivial!  
Thus, we do not know how to prove that $V_{u_0}$ has non-zero
image in $\Sha(\mathscr{E}_{u_0})[2]$ for infinitely many 
of the points $u_0$ as in Theorem \ref{phipthm},
and hence we do not know if ${\rm{rank}}(\mathscr{E}_{u_0}(
\kappa_0(T))) < 3$ (let alone if
$\mathscr{E}_{u_0}(\kappa_0(T))$ has rank 1)
for infinitely many $u_0$. 

Here is the successful strategy for using
arithmetic information from the $\mathscr{E}_{u_0}$'s
to bound the dimension in (\ref{eta2}).
We are aiming to prove that $\mathscr{E}_{\eta}(\kappa(u,T))$ has
rank 1, and in (\ref{QdefT}) we have already found a point
of infinite order, so it suffices to bound the rank from above by 1
after replacing $\kappa$ with a finite extension $\kappa'$ that may
depend on the parameters $c, d \in \kappa^{\times}$ that were used in the
definition of $\mathscr{E}_{\eta}$.
Since $\mathscr{E}_{\eta}(\overline{\kappa}(u,T))$ is 
{\em finitely generated} (Lang--N\'eron),
we may replace $\kappa$ with a suitable finite extension
(depending on $c$ and $d$)
to reduce to the case when
$\mathscr{E}_{\eta}(\overline{\kappa}(u,T)) = 
\mathscr{E}_{\eta}(\kappa(u,T))$.
Now consider the commutative diagram of natural maps
\begin{equation}\label{specdig}
\xymatrix{
{\mathscr{E}_{\eta}(\kappa(u,T))/2\cdot
\mathscr{E}_{\eta}(\kappa(u,T))} \ar[r] \ar[d] &
{\mathscr{E}_{\eta}(\overline{\kappa}(u,T))/2\cdot
\mathscr{E}_{\eta}(\overline{\kappa}(u,T))} \ar[d]\\
{\mathscr{E}_{u_0}(\kappa_0(T))/2\cdot
\mathscr{E}_{u_0}(\kappa_0(T))} \ar[r] &
{\mathscr{E}_{\overline{u}_0}(\overline{\kappa}(T))/
2\cdot\mathscr{E}_{\overline{u}_0}(\overline{\kappa}(T))}}
\end{equation}
in which 
$\overline{u}_0 \in {\mathbf{A}}^1_{\kappa}(\overline{\kappa})$ 
is a choice of geometric point 
over a closed point $u_0 \in {\mathbf{A}}^1_{\kappa}$, 
and both vertical maps are defined by the valuative criterion
for properness.  Since we adjusted $\kappa$ so 
that 
$\mathscr{E}_{\eta}(\overline{\kappa}(u,T)) = 
\mathscr{E}_{\eta}(\kappa(u,T))$, 
the top side of (\ref{specdig}) is an {\em isomorphism}.
Therefore (\ref{eta2}) is at most 2 if 
\begin{itemize}
\item the right side of (\ref{specdig}) is injective for
{\em all but finitely many}
$\overline{\kappa}$-points $\overline{u}_0 \in 
{\mathbf{A}}^1_{\kappa}(\overline{\kappa})$,
\item the
image of the map along the bottom side of (\ref{specdig}) is
at most 2-dimensional for
{\em infinitely many} closed points $u_0 \in
{\mathbf{A}}^1_{\kappa}$
(equipped with one of the finitely many
choices of $\overline{\kappa}$-point
$\overline{u}_0$ over $u_0$). 
\end{itemize} 
We consider these two respective assertions as
``geometric'' and ``arithmetic'' in nature.

\begin{remark} We do not know {\em a priori} that
the left side of (\ref{specdig}) is injective for
all but finitely many (or even infinitely many) $u_0$,
though this injectivity does follow {\em a posteriori} from our proof that
$\mathscr{E}_{\eta}(F(T))$ has rank 1; the {\em a priori} 
difficulty is due to the fact that $\kappa$ is
not separably closed (see Theorem \ref{geomreduction}).  However, 
even if we did know such injectivity, it would be useless because
our rank bounds 
for $\mathscr{E}_{u_0}(\kappa_0(T))$ are not good enough. 
The purpose of
considering (\ref{specdig}) is precisely to circumvent
our lack of information concerning the groups
$\mathscr{E}_{u_0}(\kappa_0(T))$.
\end{remark}

We shall now undertake the geometric part of the argument (injectivity 
of the right side of (\ref{specdig}) for all but finitely 
many $\overline{u}_0$).
This will be deduced 
from a more general specialization result for abelian varieties.
Let us  isolate the essential geometric properties
of $\mathscr{E}_{\eta}$ before we pass to an axiomatized
setup with an abelian variety.
Consider the surface $S = \PP^1_{\kappa} \times \PP^1_{\kappa}$ with 
factors having respective
coordinates $u$ and $T$.  
By general ``smearing out'' principles, 
$\mathscr{E}_{\eta}$ extends to an elliptic curve
$\mathscr{E}_V$ over a dense open $V \subseteq S$.  (In fact, there
is a unique maximal such open $V$, containing all others, 
and the elliptic curve
$\mathscr{E}_V$ extending $\mathscr{E}_{\eta}$ over this $V$ is unique.
This follows from a general lemma of
Faltings \cite[\S2,~Lemma~1]{faltings},
but we do not need it.)  Pick some choice of $V$
and $\mathscr{E}_V$.  There are
finitely many (if any) codimension-1 points
in $S$ not in $V$, and if $\mathscr{E}_{\eta}$
has good reduction at such a point $s$ then we can ``smear out''
the proper N\'eron model over $\OO_{S,s}$ 
and glue it to $\mathscr{E}_V$ so as to increase
$V$ to contain $s$.  Doing this finitely many times,
we may assume $V$ contains all codimension-1 points
of $S$ where $\mathscr{E}_{\eta}$ has good reduction.

The complement $S - V$ consists of finitely
many curves and isolated closed points.
Since $\mathscr{E}_{u_0}$ is smooth for all closed
points $u_0 \in {\mathbf{A}}^1_{\kappa}$, 
the curves in the complementary locus 
$$S - V \subseteq \PP^1_{\kappa} \times \PP^1_{\kappa}$$
are ``non-vertical'' except for possibly 
$\{\infty\} \times \PP^1_{\kappa}$.
Put in geometric terms, when 
the bad locus for
$\mathscr{E}_{\eta}$ over $S$
is fibered over the $T$-line it ``moves'' in
the fibers $S_t = \PP^1$ except for possibly
at the point $\infty$ in these fibers.
We need to analyze the situation along the vertical line
$u = \infty$.

\begin{lemma}\label{baddiv} The elliptic curve $\mathscr{E}_{\eta}$ 
in $(\ref{Wut})$
has bad reduction at the codimension-$1$ 
generic point $\eta_{\infty}$ of the line $u = \infty$ in $S$, 
with reduction type that is potentially good.
The ramification of $\mathscr{E}_{\eta}[2]$ at
$\eta_{\infty}$ is tame.
\end{lemma}

\begin{proof}
Since $\deg_u(h(T^2+u)) = 2p$, 
we see from (\ref{Dc4wut}) 
that $\deg_u(\Delta) = 18p$ is not divisible by 12. 
Therefore, there is bad reduction at
$\eta_{\infty}$.  The $j$-invariant $j(\mathscr{E}_{\eta})$
is a unit at $\eta_{\infty}$ because
$j$ in (\ref{Delta}) is a unit at $\infty$, so 
the reduction at $\eta_{\infty}$
is potentially good.  Since $\mathscr{E}_{\eta}[2](\eta_{\infty}) \ne O$
and the residue characteristic at $\eta_{\infty}$ is not 2, the
2-torsion $\mathscr{E}_{\eta}[2]$ is tamely ramified
at $\eta_{\infty}$.
\end{proof}

Now we pass to a general situation that uses
the properties proved in Lemma \ref{baddiv}.
Let $k$ be a separably closed field and 
let $S$ be a connected geometrically-normal
$k$-scheme of finite type,
equipped with a surjective $k$-morphism $S \rightarrow \PP^1_k$
whose fibers are geometrically reduced and whose
generic fiber is geometrically irreducible. 
By \cite[IV$_3$,~9.7.7]{ega} there is a dense
open in $\PP^1_k$ over which 
$S$ has geometrically integral fibers.
In the above discussion, 
$k = \overline{\kappa}$ and $S$ is the product
of the projective $u$-line and projective $T$-line
over $k$ with projection $S \rightarrow \PP^1_k$ onto the $u$-line.

Let $A$ be an abelian variety of
dimension $g \ge 1$ over the function field
$k(S)$.
For all but finitely many closed points $u \in \PP^1_k$,
$A$ has good reduction $A_{\eta_{u}}$ at the codimension-1
generic point $\eta_{u}$ of the geometrically integral
fiber $S_{u}$ in
the normal $S$; we write $k(S_{u})$ to denote 
the function field of this fiber.  By
the valuative criterion for properness we have a
specialization mapping
$$\rho_{u}:A(k(S)) \rightarrow A_{\eta_{u}}(k(S_{u}))$$
for such $u$.  (Since we are not assuming
that the Chow $k(S)/k$-trace of $A$ vanishes, $A(k(S))$ might not be
finitely generated.  Hence, $\rho_{u}$ cannot be defined
by elementary denominator-chasing with a finite set of
elements and their relations in $A(k(S))$, so we really do need
the valuative criterion for properness in order to define 
$\rho_{u_0}$; more specifically we cannot expect $A(k(S))$ 
to ``smear out'' beyond the codimension-1 local ring on $S$ at the generic
point $\eta_u$ of $S_u$.) 
Motivated by the goal of proving that
the right side of (\ref{specdig})
is injective with only finitely many exceptions,
we want to analyze the kernel of the reduced map
$$\rho_{u} \bmod n:   
\,\,\,A(k(S))/n\cdot A(k(S)) \rightarrow
A_{\eta_{u}}(k(S_{u}))/n\cdot A_{\eta_{u}}(k(S_{u}))$$
for suitable integers $n$ and for $u$ avoiding
a finite set of closed points on $\PP^1_k$.
To this end, it is convenient to first prove a general
finiteness lemma.

\begin{lemma}\label{chow} Let $V$ be a geometrically integral
variety over a 
field $k$ and let $B$ be an abelian variety over $K = k(V)$.
For all non-zero integers $m$ with
$\chr(k) \nmid m$, the group $B(K)/m\cdot B(K)$ is finite
if $k$ is separably closed. 
The same holds for arbitrary non-zero integers $m$ if $k$ is algebraically
closed.
\end{lemma}

\begin{proof}
We shall use Chow's theory of the $K/k$-trace
\cite[Ch.~VIII]{langav}.  Here are the key points of this theory
(for our purposes). 
In the category of pairs $(B_0,f_0)$ consisting of
an abelian variety $B_0$ over $k$ and a map $f_0:(B_0)_K \rightarrow B$
of abelian varieties over $K$, there is a final object
$({\rm{Tr}}_{K/k}(B), \tau)$ and the canonical map
$\tau:({\rm{Tr}}_{K/k}(B))_K \rightarrow B$
has infinitesimal kernel.   This object is the
{\em $K/k$-trace} of $B$.  Obviously the map
$${\rm{Tr}}_{K/k}(B)(k) \hookrightarrow
{\rm{Tr}}_{K/k}(B)(K) \stackrel{\tau}{\rightarrow} B(K)$$
is injective.  The
Lang--N\'eron theorem \cite[Thm.~1]{langneron} says that the quotient group
\begin{equation}\label{bquot}
B(K)/{\rm{Tr}}_{K/k}(B)(k)
\end{equation}
is finitely generated.  
(To the best of our knowledge, all published references on
these topics are written in pre-Grothendieck terminology; 
the reader is referred to \cite{bnotes} for a discussion
of the Chow trace and Lang--N\'eron theorem
using scheme-theoretic methods.)

Now assume that $k$ is separably closed.
Since ${\rm{Tr}}_{K/k}(B)(k)$ is 
the group of rational points of an abelian variety over
a separably closed field, it is $m$-divisible
(and the restriction $\chr(k) \nmid m$
can be removed if $k$ is algebraically closed).  Thus,
$$
{\rm{Tr}}_{K/k}(B)(k) \subseteq m \cdot B(K),
$$
so 
$$B(K)/m \cdot B(K) \simeq (B(K)/{\rm{Tr}}_{K/k}(B)(k))/m
\cdot (B(K)/{\rm{Tr}}_{K/k}(B)(k)).$$
This yields the desired finiteness because 
(\ref{bquot}) is finitely generated.
\end{proof}

We return to the abelian variety $A_{/k(S)}$, described before 
Lemma \ref{chow}.

\begin{theorem}\label{geomreduction}
Assume that $k$ is separably closed and that for all closed points
$u \in \PP^1_k$ distinct from
$\infty$, $A$ has good reduction at some generic
point of the $($possibly reducible$)$ geometrically-reduced
fiber 
$S_{u}$.  Assume moreover that at
some generic point $\eta_{\infty}$
of the fiber $S_{\infty}$
there is potentially good reduction. 

Fix $n \in \ZZ$ with $|n| > 1$
such that $\chr(k) \nmid n$, and assume
that the Galois splitting field of
the finite \'etale $k(S)$-group $A[n]$ 
is tamely ramified at the codimension-$1$ 
point $\eta_{\infty} \in S$.

The mod-$n$ reduction
$$\rho_{u} \bmod n:\,\,\,A(k(S))/n\cdot A(k(S)) \rightarrow
A_{\eta_{u}}(k(S_{u}))/n\cdot A_{\eta_{u}}(k(S_{u}))$$
of the specialization map along $S_u$ is injective 
for all but finitely many closed points $u \in \PP^1_k$.
\end{theorem}

The tameness assumption is equivalent to the
condition that $A$ acquires good reduction
over a finite separable extension of $k(S)$
that is tame at a place over $\eta_{\eta}$
(this is explained in the proof),
and so this hypothesis is automatically satisfied
when every positive prime $\ell \le 2g+1$ is a unit in $k$
(that is, $\chr(k) = 0$ or $\chr(k) > 2g+1$).
Thus, by setting $g = 1$ and $n = 2$ in Theorem
\ref{geomreduction}, we may conclude via Lemma \ref{baddiv} 
(which also gives the desired tameness in characteristic 3) that
the right side of
(\ref{specdig}) 
is injective for all but finitely many
$\overline{u}_0 \in {\mathbf{A}}^1_{\kappa}(\overline{\kappa})$.

\begin{proof}
The hypotheses on $S$ and $A$ are preserved
under extension of the base field.  Moreover, if
$\overline{k}$ is an algebraic closure of $k$ then we
claim that the natural map
$$A(k(S))/n \cdot A(k(S)) \rightarrow A(\overline{k}(S))/n \cdot
A(\overline{k}(S))$$
is injective, so we may reduce to the case when $k$ is
algebraically closed.  The case of characteristic 0 is trivial, so
we can assume ${\rm{char}}(k) = p > 0$.   
It suffices to check more generally that if
$K$ is a field with characteristic $p > 0$ and
$G$ is a commutative $K$-group of finite type then 
the map $G(K)/n \cdot G(K) \rightarrow G(K')/n \cdot G(K')$
is injective for any purely inseparable algebraic extension $K'/K$ 
and any integer $n$ not divisible by $p$. We may
assume $K' = K^{1/p}$, so we get an identification
$G(K') \simeq G^{(p)}(K)$ that identifies
the inclusion $G(K) \rightarrow G(K')$ with the map on 
$K$-points induced by
the relative Frobenius morphism $F_G:G \rightarrow G^{(p)}$.
Since $[p]:G \rightarrow G$ factors through $F_G$
\cite[VII$_{\rm{A}}$,~\S4.3]{sga3}, it suffices
to prove that the $p$-torsion in $G(K)/n\cdot G(K)$ vanishes,
and this is clear since $p \nmid n$.

Let $W \subseteq S$ be a dense open such that 
$A$ extends to an abelian scheme
$A_W$ over $W$.  The complement
$S - W$ contains at most finitely many codimension-1
points of $S$, and if $A$ has good reduction
at any such point $s$ then we may glue $A_W$
with a smearing-out of the proper N\'eron model
of $A$ over ${\mathscr{O}}_{S,s}$ to increase $W$ to contain
a neighborhood of $s$.  Thus, by the hypothesis on reduction for $A$,
we may suppose that 
no fiber $S_{u}$ over a closed
any point $u \in \PP^1_k$ is disjoint 
from $W$ except for possibly
$S_{\infty}$.  This property
of $W$ is unaffected by shrinking $W$ in codimension $\ge 2$.
Let $\eta$ be the generic point of $S$. 

By Lemma \ref{chow} with $V = S$, $A(k(S))/n \cdot A(k(S))$ is finite.
We conclude from the pigeonhole principle 
that if $\rho_{u} \bmod n$
has nontrivial kernel for infinitely many $u$
(ignoring the finitely many for which 
$S_{u}$ is reducible, in which case
$\rho_{u}$ is not defined), then
some non-zero
$$\overline{R}
 \in A(k(S))/n\cdot A(k(S))$$ is killed by $\rho_{u} \bmod n$
for infinitely many $u$.  
Thus, it suffices to prove that
if $R_{\eta} \in A(k(S))$
has the property that $\rho_{u}(R_{\eta})$
lies in $n \cdot A_{\eta_{u}}(k(S_{u}))$
for infinitely many $u$ (ignoring
the finitely many $u$ for which
$\rho_{u}$ is not defined) then $R_{\eta} \in n \cdot A(k(S))$. 

Choose $R_{\eta} \in A(k(S))$ such that $\rho_{u}(R_{\eta})$ 
lies in $n \cdot
A_{\eta_{u}}(k(S_{u}))$ for infinitely many 
$u$. By denominator-chasing, 
$R_{\eta}$ extends (uniquely) to $R_U \in A_W(U)$ for some
dense open $U \subseteq W$.
The valuative criterion
for properness extends $R_U$
over each of the finitely many 
codimension-1 points of $W$ not contained in $U$.
Thus, by shrinking $W$ in codimension $\ge 2$ if necessary, we may
assume that $R_{\eta}$ extends to a section $R \in A_W(W)$ of
the abelian scheme $A_W \rightarrow W$. 

The pullback of $[n]:A_W \rightarrow A_W$ along $R \in A_W(W)$
is a finite \'etale cover
\begin{equation}\label{nrw}
[n]^{-1}(R) \rightarrow W.
\end{equation}
Our goal is to prove that (\ref{nrw})
has a section
over the generic point
$\eta = \Spec k(S)$ of $W$.
Let $L$ be a residue field on 
$[n]^{-1}(R)_{\eta} = [n]^{-1}(R_{\eta})$,
so $L$ is a finite separable extension of
$k(S)$, say with degree $d_L$.  We want $d_L = 1$ for some such $L$.

\begin{lemma}\label{algclosed}
For each $L$, the subfield $k(\PP^1)$ is algebraically closed in $L$.
\end{lemma}

\begin{proof} 
Let $K/k(\PP^1)$ be the algebraic closure of
$k(\PP^1)$ in $L$, so $K/k(\PP^1)$ is a finite
separable extension because $L/k(\PP^1)$ is a finitely
generated separable extension
(as $k(S)$ is separable
over $k(\PP^1)$, since the generic fiber of
$S \rightarrow \PP^1$ is geometrically integral).  
The intermediate
fields $K$ and $k(S)$ in the separable extension $L/k(\PP^1)$ are linearly
disjoint over $k(\PP^1)$ because $K/k(\PP^1)$ is algebraic
and $k(\PP^1)$ is algebraically closed in
$k(S)$.  Thus, if $\theta \in \PP^1$ is the generic
point then the function field
\begin{equation}\label{kxs}
K(S_{\theta}) \eqdef K \otimes_{k(\PP^1)} k(S)
\end{equation}
of the geometrically integral
generic fiber ${S_{\theta}}_{/K}$ is identified with the
intermediate composite field $K \cdot k(S)$ in $L/k(S)$. 
The hypothesis on the good reduction of $A$ implies
that for every closed point $u \in \PP^1 - \{\infty\}$, 
some generic point $\eta_{u}$ of the reduced fiber 
$S_{u}$ lies in $W$. 
Hence, since  $[n]^{-1}(R) \rightarrow W$ is a finite
\'etale cover, 
the residue field $L$ on 
$[n]^{-1}(R_{\eta})$
is unramified over
the discrete valuation on $k(S)$ arising from some such 
$\eta_{u}$ for every closed point $u \in \PP^1 - \{\infty\}$.
It follows that for every such $u$, the
intermediate finite separable extension $K(S_{\theta})/k(S)$
is also unramified at some generic point $\eta_{u}$ of
$S_{u}$. 

We also need to understand the ramification behavior of
$L/k(S)$ at the discrete valuation on $k(S)$ arising from
a generic point $\eta_{\infty}$ on the reduced fibral curve 
$S_{\infty}$ such that
$A$ has potentially good reduction over a tame extension
at $\eta_{\infty}$; the existence of such an $\eta_{\infty}$
was one of our initial assumptions on $A$. 
We claim that $L/k(S)$ is tamely ramified
at {\em all} places of $L$ over $\eta_{\infty}$. 
Some care will be required because
$L/k(S)$ may be non-Galois.  

The first step is to check
that $L$ admits at least one place that is tame over $\eta_{\infty}$,
and to do this it suffices to choose a separable closure
of the residue field at $\eta_{\infty}$ and to show that $[n]^{-1}(R_{\eta})$
splits over a tame extension of 
the fraction field of the
associated strict henselization $\OO_{S,\eta_{\infty}}^{\rm{sh}}$.
Since $A[n]$ is assumed to be tamely ramified at
$\eta_{\infty}$, there exists 
a finite tame extension $F'$ over
the fraction field of $\OO_{S,\eta_{\infty}}^{\rm{sh}}$
such that $A[n]_{F'}$ is a constant group.
We can assume $|n| > 1$, so there exists a prime $\ell|n$
and $\ell \ne {\rm{char}}(k)$.  Since $A[n]_{F'}$ is constant, the
Galois-action on the $\ell$-adic
Tate module of $A_{/F'}$ 
has pro-$\ell$ image that is finite
(since $A$ has potentially good reduction at $\eta_{\infty}$),
so after replacing $F'$ with a suitable $\ell$-power
(hence {\em tame}) extension
we can assume that $A_{/F'}$ has good reduction.
Let $\mathscr{A}$ denote
the proper N\'eron model of $A_{/F'}$ over
the integral closure $\OO_{F'}$ of
$\OO_{S,\eta_{\infty}}^{\rm{sh}}$ in $F'$.
The group $A(F') = \mathscr{A}(\OO_{F'})$ is $n$-divisible because 
$\OO_{F'}$ is strictly henselian and $n$ is
not divisible by the residue characteristic of $\OO_{F'}$, so 
$[n]^{-1}(R_{\eta})(F') \ne \emptyset$. Since $A[n]_{F'}$
is split, it follows that the \'etale $A[n]$-torsor
$[n]^{-1}(R_{\eta})$ must therefore be split over $F'$.
Hence, $L$ admits a $k(S)$-embedding into $F'$,
so 
$L/k(S)$ is tamely ramified at some place
$w_L$ over 
the discrete valuation on $k(S)$ arising from 
$\eta_{\infty}$.

By definition, 
$L$ is a residue field on 
an \'etale $A[n]$-torsor
$[n]^{-1}(R_{\eta})$ over $k(S)$,
and (by hypothesis) the $k(S)$-group $A[n]$ splits
over a finite Galois extension $M/k(S)$ that is
tamely ramified over ${\eta}_{\infty}$.  Thus, 
the factor fields of the finite \'etale $M$-algebra $L \otimes_{k(S)} M$ are
residue fields on the torsor $[n]^{-1}(R_{\eta})_M$
for a finite constant group
over $M$ (namely, the constant
group $A[n]_M$).  Hence, the factor
fields $L_i$ of $L \otimes_{k(S)} M$ are Galois over $M$
and the $L_i$'s are pairwise $M$-isomorphic.  Pick a 
place $w_M$ on $M$ lifting the place $\eta_{\infty}$ on $k(S)$.
Since $w_M$ and $w_L$ lift the same place on $k(S)$, 
we can find 
a factor field $L_{i_w}$ of $L \otimes_{k(S)} M$ and a place $v_{i_w}$ on
$L_{i_w}$ that lifts the places $w_L$ 
and $w_M$.  The place $v_{i_w}$ on $L_{i_w}$ must be tame
over the place $w_M$ because $w_L$ 
is tame over $\eta_{\infty}$ on $k(S)$.  The extension
$L_{i_w}/M$ is Galois, so $L_{i_w}/M$ is tame
at {\em all} places over $w_M$.   
Since the $L_i$'s are pairwise $M$-isomorphic and $w_M$ is 
an arbitrary 
place on $M$ over $\eta_{\infty}$, every $L_i$ is tame over every 
place on $M$ lifting 
$\eta_{\infty}$.
Since {\em all} places
of $M$ over $\eta_{\infty}$ are tame over $\eta_{\infty}$,
we conclude that all places lying over
$\eta_{\infty}$ on each $L_i$
are tame over $\eta_{\infty}$.  Upon choosing some
$L_{i_0}$, 
the extension $L/k(S)$ is a subextension of 
$L_{i_0}/k(S)$ and hence 
$L/k(S)$ is tamely ramified at {\em all} places
over $\eta_{\infty}$. 
The same therefore holds for the intermediate extension
$K(S_{\theta})/k(S)$ in (\ref{kxs}). 

Summarizing our conclusions, 
the finite separable extension $K(S_{\theta}) = K \otimes_{k(\PP^1)} k(S)$
over $k(S)$ is unramified at some generic point
of the reduced fiber $S_{u}$ for each $u \ne \infty$
and is tamely ramified over some generic point of
the reduced fiber $S_{\infty}$.  The reducedness of
the fibers implies that a uniformizer at
a closed point $u \in \PP^1$ pulls back to be
a uniformizer in the local ring at the codimension-1
point $\eta_{u}$ on the normal surface $S$. 
Hence, the discrete valuation on $k(\PP^1)$
associated to $u$ has ramification index 1
under the discrete valuation on $k(S)$
associated to $\eta_{u}$, and the corresponding
residue field extension is separable (because the
residue field at $u$ is the field $k$ that is
algebraically closed).  It follows by classical
valuation theory and (\ref{kxs}) that
if $\eta_{u}$ is unramified (resp. tamely ramified)
in $K(S_{\theta})$ then $u$ must be unramified
(resp. tamely ramified) in $K$.  Hence, the
finite separable (possibly non-Galois) extension 
$K/k(\PP^1)$ is unramified away from $\infty$
and is tamely ramified at {\em all} places over
$\infty$.  Since $k$ is separably closed,
we conclude that $K = k(\PP^1)$.
\end{proof}

We return to the proof of Theorem \ref{geomreduction}.
Let  $\mathscr{C}_L$ be the connected component of 
$[n]^{-1}(R)$ with function field $L$. 
Since $L/k(\PP^1)$ is a finitely generated
separable extension with transcendence degree 1,
it follows from Lemma \ref{algclosed}
that the fiber of $\mathscr{C}_L$ over
the generic point of $\PP^1_k$ must be geometrically integral 
over $k(\PP^1)$. 
Hence, by
\cite[IV$_3$,~9.7.7]{ega}, 
there is a Zariski-dense open $U_L \subseteq
{\mathbf{A}}^1_k$ such that the fiber $(\mathscr{C}_L)_{u}$ is
geometrically integral over $k(u)$ for all $u \in U_L$. 
By removing finitely many closed points from $U_L$, we may
(and do) also assume that $S_{u}$ is geometrically
integral over $k(u)$ for all $u \in U_L$. 
Since $[n]^{-1}(R)$ is finite \'etale over $W$ and the open subset 
$W_{u} \subseteq S_{u}$ is non-empty
for all $u \in {\mathbf{A}}^1_k$, 
the finite \'etale map 
$(\mathscr{C}_L)_{u} \rightarrow
W_{u}$ has degree $$[k(\mathscr{C}_L):k(W)] = [L:k(S)] = d_L$$
for all points $u \in U_L$.

Choose
$u \in \cap_L U_L$, 
where $L$ runs over all the residue fields on $[n]^{-1}(R_{\eta})$.
We have just proved that
the fiber $(\mathscr{C}_L)_{u}$ is connected 
(even geometrically integral over $k(u)$) for all $L$. 
It follows that $\{(\mathscr{C}_L)_{u}\}_L$ is
the set of connected components of the finite \'etale $W_{u}$-scheme 
$[n]^{-1}(R)_{u} = [n]^{-1}(R_{u})$ and the map 
$$(\mathscr{C}_L)_{u} \rightarrow W_{u} \subseteq S_{u}$$
is \'etale with
generic degree $d_L$ for all $L$.
Membership in $\cap_L U_L$ omits only finitely many closed points 
$u$, so by the hypothesis 
$\rho_{u}(R_{\eta}) \in n \cdot A_{\eta_{u}}(k(S_{u}))$
for infinitely many $u$ (with $\rho_{u}(R_{\eta})$
the generic point of the section $R_{u}$ of $A_W$ over $W_{u}$)
we conclude that there exists
a closed point $u' \in \cap_L U_L$
such that $$\rho_{u'}(R_{\eta}) \in n \cdot A_{\eta_{u'}}(k(S_{u'})).$$
In particular, the $W_{u'}$-\'etale scheme 
$[n]^{-1}(R_{u'})$ has a $k(S_{u'})$-rational point. 
This rational point lies in some
fibral connected component $(\mathscr{C}_{L_0})_{u'}$, so 
the generic degree $d_{L_0}$ of this component
over $S_{u'}$ must equal 1.
\end{proof}

\section{Generic rank bound II. Arithmetic arguments}\label{dull}

Our remaining task is to prove
that the bottom side of (\ref{specdig}) is
injective for infinitely many closed points
$u_0 \in {\mathbf{A}}^1_{\kappa}$.
In Theorem 
\ref{phipthm} we will find 
infinitely many points $u_0$ such that the elliptic curve 
${\mathscr E}_{u_0/\kappa_0(T)}$ 
over the global field
$\kappa_0(T)$ 
has exactly two places of bad reduction,
and in \S\ref{selmer} we will 
prove injectivity along the bottom of (\ref{specdig}) 
for such $u_0$. 

As preparation for the study of the image along the bottom side of
(\ref{specdig}) for well-chosen closed
points $u_0 \in {\mathbf{A}}^1_{\kappa}$,
we fix an arbitrary $u_0$ and
find the reduction
type of $\mathscr{E}_{u_0}$ at each place of 
$\kappa_0(T)$.  After we find these
reduction types, the points $u_0$ that will become our
focus of interest will be those
such that ${\mathscr E}_{u_0}$ has the smallest possible number of 
physical points
of bad reduction on the $T$-line
$\PP^1_{\kappa_0}$.

Recall that (\ref{Wut}) 
defines $\mathscr{E}_{\eta}$ in terms of $h(T^2+u)$, where 
$h(T) = c T^{2p} + du$.  From (\ref{Dc4wut}), the discriminant of 
(\ref{Wut}) involves 
$h(T^2+u)$ and $1+4h(T^2+u)$. 
Under a $u_0$-specialization,
$h(T^2+u)$ becomes a $p$th power in
$\kappa_0[T]$:
$$
h(T^2+u)|_{u=u_0} =
(c^{1/p} (T^2 + u_0)^2 + d^{1/p}u_0^{1/p})^p.
$$
Likewise, $1 + 4h(T^2+u)$ specializes to a $p$th power in $\kappa_0[T]$:
$$
(1+4h(T^2+u))|_{u=u_0} = (1 + 4(c^{1/p} (T^2+u_0)^2 + d^{1/p}u_0^{1/p}))^p.
$$
For all but finitely many closed
points $u_0 \in {\mathbf{A}}^1_{\kappa}$, the 
$p$th-root polynomials
\begin{equation}\label{pijsdef}
\pi_1 := 
c^{1/p} (T^2+u_0)^2 + d^{1/p}u_0^{1/p},\,\,\,
\pi_2 := 1 + 4\pi_1
\end{equation}
are separable in $\kappa_0[T]$.  (These quartics over the finite
field $\kappa_0$ 
may be reducible for many points $u_0$, and so
even in characteristic $p > 3$ these quartics may be fail
to be separable for some non-empty finite set of points $u_0$.
In Theorem \ref{phipthm} below, we will show that for infinitely
many $u_0$ we can do much better than mere separability.)
We now restrict attention to
those $u_0$ such that the two polynomials in (\ref{pijsdef}) 
are both separable. 
(Our notation $\pi_1$ and $\pi_2$ does not indicate 
the dependence on $u_0$; it would be more accurate to write 
$\pi_{1,u_0}$ and $\pi_{2,u_0}$, but we simply ask the reader to 
remember the
dependence on $u_0$.)

Specializing (\ref{Dc4wut}) at $u_0$, the Weierstrass model
that defines ${\mathscr E}_{u_0/\kappa_0(T)}$
 has parameters 
\begin{equation}\label{delc4}
\Delta|_{u=u_0} = 16\pi_1^{8p}\pi_2^{p}, \ \ \ 
c_4|_{u=u_0} = 16\pi_1^{2p}(1+3\pi_1^p) = 
16\pi_1^{2p}(\pi_2 -\pi_1)^p,  
\end{equation}
and this Weierstrass model is
integral away from $T = \infty$.
Thus, the only possible bad reduction for $\mathscr{E}_{u_0}$ 
over the $T$-line $\PP^1_{\kappa_0}$ is 
at $\infty$ and 
at the zeros of $\pi_1$ and $\pi_2$.

What is the behavior of
$\mathscr{E}_{u_0}$ at the point $\infty \in \PP^1_{\kappa_0}$? 
We return to Lemmas \ref{good} and \ref{bad}. 
Both $\pi_1$ and $\pi_2$ have degree $4$ in $\kappa_0[T]$, so 
by (\ref{delc4}) we have 
$$
\ord_\infty(\Delta|_{u=u_0}) = -36p.
$$
When $\chr(\kappa) > 3$, 
$$
\ord_\infty(c_4|_{u=u_0}) = -12p, \ \
\ord_\infty(j|_{u=u_0}) = 0.
$$
When $\chr(\kappa) = 3$, 
$$
\ord_\infty(c_4|_{u=u_0}) = -8p, \ \
\ord_\infty(j|_{u=u_0}) = 12p.
$$
Thus, there is potentially good reduction at $T=\infty$
in all cases, and
Lemma \ref{good} ensures that this reduction is good.

Now we analyze the reduction types at points $x_j$ in
the zero-scheme of $\pi_j$ on $\PP^1_{\kappa_0}$. 
Since $\pi_1$ is separable in $\kappa_0[T]$ (by our choice of $u_0$), 
we see from (\ref{delc4}) that 
for any $x_1$,  
$$
\ord_{x_1}(\Delta|_{u=u_0}) = 8p, \ \ \
\ord_{x_1}(c_4|_{u=u_0}) = 2p \equiv 2 \bmod 4.
$$
Therefore $\ord_{x_1}(j(\mathscr{E}_{u_0})) = 6p-8p=-2p < 0$. 
By Lemma \ref{bad}, there must be
(potentially multiplicative) additive reduction
at $x_1$.
Similarly, we compute 
$$
\ord_{x_2}(\Delta|_{u=u_0}) = p, \ \ \
\ord_{x_2}(c_4|_{u=u_0}) = 0,
$$
so $\ord_{x_2}(j(\mathscr{E}_{u_0})) = -p < 0$.  By Lemma \ref{bad},
there is multiplicative reduction at $x_2$.

We have shown that
the N\'eron model $N(\mathscr{E}_{u_0}) \rightarrow
\PP^1_{\kappa_0}$ enjoys the following reduction properties:
\begin{list}
{({{\alph{bean}}})}
{\usecounter{bean}
\setlength{\rightmargin}{\leftmargin}}
\item good reduction at
all closed points of $\PP^1_{\kappa_0}$ away from 
zeros of $\pi_1$ and $\pi_2$,
\item multiplicative reduction
at zeros $x_2$ of $\pi_2$, with $\ord_{x_2}(j_{u_0}) = -p$. 
\item additive reduction at zeros $x_1$ of $\pi_1$, with 
$\ord_{x_1}(j_{u_0}) = -2p < 0$.
\end{list}
Properties (b) and (c) will be used in our work with
N\'eron models and Selmer groups in \S\ref{selmer}, but
now we focus on (a).  
The most favorable $u_0$'s for our purposes will be those
such that $\mathscr{E}_{u_0}$ has the least
possible number of physical points of bad reduction,
so we want to find many $u_0$ such that 
$\pi_1$ and $\pi_2$ are both irreducible in 
$\kappa_0[T]$.  For
such $u_0$,  
$\mathscr{E}_{u_0}$ has exactly two physical points
of bad reduction on $\PP^1_{\kappa_0}$.

\begin{theorem}\label{phipthm}
There exist infinitely many closed points $u_0 \in
{\mathbf{A}}^1_{\kappa}$ such
that $\pi_1, \pi_2 \in \kappa_0[T]$ are irreducible.
\end{theorem}

To find the infinitely many $u_0$
as in the theorem will require some effort,
so let us first sketch the basic idea.
In (\ref{pijsdef}) we see that $u_0 \in \kappa_0$ intervenes in $\pi_1$ and $\pi_2$
through the value $u_0^{1/p} \in \kappa_0$, so to put ourselves
in the position of specializing polynomials in $u$
we apply arithmetic Frobenius of $\kappa_0$
to the coefficients of $\pi_1$ and $\pi_2$.  This leads
us to consider the polynomials
\begin{equation}\label{bigpi}
\Pi_1(u,T) := c(T^2 + u^p)^2 + d u, \,\,\, \Pi_2(u,T) := 1 + 4\Pi_1(u,T)
\in \kappa[u][T].
\end{equation}
For any closed point $u_0 \in {\mathbf{A}}^1_{\kappa}$, 
the specialization $\Pi_j(u_0,T) \in \kappa_0[T]$ is
the image of $\pi_j \in \kappa_0[T]$ under the arithmetic
Frobenius automorphism of $\kappa_0$.  Thus, 
Theorem \ref{phipthm} is 
equivalent to the existence of infinitely many
$u_0 \in {\mathbf{A}}^1_{\kappa}$ 
such that both $\Pi_1(u_0,T)$ and $\Pi_2(u_0,T)$
are irreducible in $\kappa_0[T]$, where
$\kappa_0 = \kappa(u_0)$ is varying with $u_0$. 
It
is this equivalent statement that we will 
actually prove (Theorem \ref{newthm} below).

Expanding $\Pi_1$ and $\Pi_2$ as polynomials in $\kappa(u)[T]$, 
we have
\begin{equation}\label{Pi12}
\Pi_1 = c \left(T^4 + 2 u^p T^2 + u^{2p} 
+ \frac{du}{c}\right),\,\,\,
\Pi_2 = 4c \left(T^4 + 2 u^p T^2 + u^{2p} + \frac{du}{c} +
\frac{1}{4c}\right).
\end{equation}
It is left to the reader to check that $\Pi_1$ and $\Pi_2$
are separable and irreducible over $\kappa(u)$, via the following
elementary criterion concerning polynomials of the form $X^4 + aX^2 + b$.

\begin{lemma}\label{irred4}
Let $K$ be a field with ${\rm{char}}(K) \ne 2$.
A polynomial $f = X^4 + aX^2 + b \in K[X]$ 
is
separable if and only if $b$ and $a^2 - 4b$ are non-zero. 
It is irreducible
if $b$ and $a^2 - 4b$
are non-squares in $K^\times$.
\end{lemma}

\begin{proof} 
The condition for separability is obvious.  We now assume that $b$ and
$a^2 - 4b$ are non-squares in $K^{\times}$.  Since $a^2 - 4b$
is not a square, $f$ has no roots in $K$ and has
no factors of the form $X^2 - c$ in $K[X]$.
Thus, if we write the four
roots of $f$ in a splitting field
as $\pm r_1$ and $\pm r_2$, a non-trivial monic
factor of $f$ in $K[X]$ must have the form 
$(X \pm r_1)(X \pm r_2)$. If such a factor exists
then $r_1r_2 \in K$ and 
$b = (r_1r_2)^2$, contradicting the assumption that 
$b$ is a non-square in $K$.  
\end{proof}

In view of the irreducibility of 
each $\Pi_j$ in $\kappa(u)[T]$ and our desire
to prove $$\Pi_1(u_0,T), \,\, \Pi_2(u_0,T) \in \kappa_0[T]$$
are irreducible for infinitely many closed points $u_0 \in
\mathbf{A}^1_{\kappa}$, our problem 
resembles Hilbert irreducibility.  However, 
finite fields are not Hilbertian and anyway we
are not generally specializing $u$ at elements of $\kappa$
(since $[\kappa_0:\kappa] > 1$ with only finitely
many exceptions).

The main idea that will produce the desired $u_0$'s is 
the following theorem.  It gives
a group-theoretic criterion for a polynomial over
a global field to specialize to an irreducible polynomial
over the residue field at infinitely many places
(see Remark \ref{chebrem}).

\begin{theorem}\label{keithfrobold}
Let $K$ be a global field and let
$f \in K[T]$ be a monic separable irreducible
polynomial of degree $n$.  Let $K'/K$ be a splitting field for $f$
and let $G = {\rm{Gal}}(K'/K)$.
For any non-archimedean place $v$ of $K$ at which
$f$ has integral coefficients, 
$f \bmod v$ is irreducible over the residue field ${\mathbf F}_v$ at $v$
if and only if $v$ is unramified in $K'$ and the Frobenius elements over 
$v$ 
in $G$ act as $n$-cycles on the set of roots of $f$ in $K'$.
\end{theorem}

\begin{proof}
Let $r$ be a root of $f$ in $K'$. 
If $f$ is $v$-integral and $f \bmod v$ is separable, then the
discriminant of $f$ is a $v$-adic unit, so $v$ is unramified 
in $K(r)$.
Since $K'$ is a composite of such extensions of $K$, 
in such cases $v$ must be unramified in 
$K'$. Let $v'$ be a place of $K'$ over a
place $v$ in $K$ that is unramified in $K'$.
The action of ${\rm{Frob}}(v'|v)$
on the $n$ roots of $f$ in $K'$ is identified
with the action of 
the finite-field Frobenius $x \mapsto x^{\#\FF_v}$ 
on the full set of $n$ roots of $f \bmod v$ (in $\mathbf{F}_{v'}$). 
In particular, $f \bmod v$ is irreducible 
over $\FF_v$ if and only if 
$v$ is unramified in 
$K'$ and ${\rm{Frob}}(v'|v)$ acts as an $n$-cycle on the roots of $f$.
\end{proof}

\begin{remark}\label{chebrem}
In the setting of Theorem \ref{keithfrobold}, if $r \in K'$ is
a root of $f$ and $H \subseteq G$ is the subgroup associated
to the intermediate field $K(r) \subseteq K'$, then an element
$\gamma \in G$ acts as an $n$-cycle on the set of
roots of $f$ in $K'$ if and only if the cyclic
subgroup $\langle \gamma \rangle$ is a set of
representatives for the coset space $G/H$ of order $n$.
We conclude by Chebotarev's density theorem
that 
$f \bmod v$ is irreducible for infinitely many
places $v$ of $K$ if and only if 
$G/H$ admits a set of representatives
that is a cyclic subgroup of $G$.
\end{remark}

\begin{corollary}\label{keithfrob}
Let $K$ be a global field and let $f \in K[T]$ be a monic 
separable irreducible polynomial
of degree $n$.  The following are equivalent
$($restricting attention to non-archimedean places at which
the coefficients of $f$ are integral$)$:

\begin{enumerate}
\item There is some place $v$ such that $f \bmod v$ is irreducible.

\item There is a positive proportion of places $v$ such that $f \bmod v$ is 
irreducible.
\end{enumerate}
\end{corollary}

\begin{proof}
The implication $(2) \Rightarrow (1)$ is trivial, and the converse
follows from Theorem \ref{keithfrobold} and Chebotarev's density theorem.
\end{proof}

\begin{example}\label{specex}
Let  $f$ satisfy the hypotheses in Theorem \ref{keithfrobold},
and let $\{r_1,\dots,r_n\}$ be an ordering of the set
of roots of $f$ in $K'$.  Identify $G = {\rm{Gal}}(K'/K)$ with a subgroup
$\overline{G} \subseteq S_n$ via the $G$-action on the $r_j$'s.
By Theorem \ref{keithfrobold},
$f \bmod v$ is irreducible for infinitely many
$v$ if and only if $\overline{G}$ contains an $n$-cycle. 
\begin{enumerate}

\item If $G$ is 
isomorphic to $S_n$ as abstract groups (where $n = \deg f$), 
then $\overline{G} =
S_n$.  Since $\overline{G}$ contains an $n$-cycle, 
$f \bmod v$ is irreducible for infinitely many $v$.

\item 
What if $G$ is isomorphic to $A_n$ (as abstract groups)? 
Since $A_n$ embeds into $S_n$ with only one possible image, 
and $A_n$ contains an $n$-cycle if and only if $n$ is odd, we  
see that $f \bmod v$ is irreducible 
infinitely often if and only if $n$ is odd.

\item 
What if $G$ (and thus $\overline{G}$)
is isomorphic to $D_n$ (as abstract groups) with $n > 2$? 
Then $\overline{G}$ is isomorphic to $D_n$ as a permutation group, so 
$f \bmod v$ is irreducible infinitely often.

The identification of $\overline{G}$ with $D_n$ as a 
permutation group 
was explained to us by D. Pollack.  Write $\overline{G} = 
\langle \sigma, \tau\rangle$, where $\sigma^n = 1$, $\tau^2 = 1$ and 
$\tau\sigma\tau^{-1} = \sigma^{-1}$.  
Since $\langle \sigma \rangle$ is normal in $\overline{G}$
and $\overline{G}$ is a transitive subgroup of $S_n$, 
all $\langle \sigma \rangle$-orbits have the same length.  
Therefore, since $\sigma$ has order $n$ it must be an $n$-cycle. 
Writing $\sigma = (1,2,\dots,n)$, 
the condition $\tau\sigma \tau^{-1} = \sigma^{-1}$ says 
$(\tau(1),\tau(2),\dots,\tau(n)) = (n,n-1,\dots,1)$ as $n$-cycles. 
We can replace $\tau$ in the 
presentation of $\overline{G}$ with $\tau\sigma^k$ for any $k$, 
so we may assume $\tau(1) = 1$.  Identifying 
$j$ with $e^{2\pi{i}(j-1)/n}$, $\sigma$ and 
$\tau$ are now the standard generators for $D_n$ in its natural 
action on an $n$-gon.

\item What if $f$ is a normal polynomial; \ie, $G$ has order $n$?  A 
transitive subgroup of order $n$ in $S_n$ contains an $n$-cycle 
if and only if it is cyclic, so 
the reduction of $f$ at infinitely many places is irreducible if $G$ is 
cyclic but not otherwise.

\item  
In the preceding four examples, the structure of the Galois group $G$ 
as an abstract group was sufficient to determine if 
the permutation group $\overline{G}$
contains an $n$-cycle.  However, this is not generally the case. 
For example, there is a group of size $2592 = 2^5 \cdot 3^4$ 
admitting two transitive actions of degree 12 such that
one action contains 12-cycles and the other does not. 
The actions were found for us by N. Boston using MAGMA; they 
are the 245th and 246th transitive groups of degree 12 
in MAGMA's enumeration.  MAGMA also  
realizes both of these transitive groups as Galois groups over $\QQ$.

\end{enumerate}
\end{example}

We now apply these ideas to the polynomials $\Pi_1(u,T)$ and 
$\Pi_2(u,T)$ from (\ref{Pi12}).
To determine their Galois groups over $\kappa(u)$, we
use the following classical lemma.

\begin{lemma}\label{langex}
Let $K$ be a field with $\chr(K) \ne 2$, and let
$f = X^4 + aX^2 + b \in K[X]$ be separable and irreducible.
Let $K'/K$ be a splitting field
and $G = {\rm{Gal}}(K'/K)$.  We have the following possibilities
for $G$ as an abstract group:
\begin{itemize}
\item $G \simeq \ZZ/4\ZZ$ if and only if
$b(a^2 - 4b) \in K^{\times}$ is a square, 
in which case the quadratic subfield is $K(\sqrt{b}) =
K(\sqrt{a^2 - 4b})$, 
\item $G \simeq
\ZZ/2\ZZ \times \ZZ/2\ZZ$
if and only if $b \in K^{\times}$ is a square, 
in which case the quadratic subfields are
$K(\sqrt{a^2 - 4b})$ , $K(\sqrt{-a + 2\sqrt{b}})$, and 
$K(\sqrt{-a - 2\sqrt{b}})$
for a fixed choice of $\sqrt{b} \in K^{\times}$,
\item $G \simeq D_4$ if and only if $b$ and $b(a^2-4b)$ 
are not squares in $K^\times$, 
in which case 
the quadratic subfields are
$K(\sqrt{a^2 - 4b})$, $K(\sqrt{b})$, and $K(\sqrt{b(a^2 - 4b)})$.
The unique quadratic subfield
over which $K'$ is a cyclic extension is
$K(\sqrt{b(a^2 - 4b)})$.
\end{itemize}
\end{lemma}

\begin{proof}
This classification of Galois groups according to 
properties of the coefficients 
can be found as an exercise in many basic algebra books, 
although usually it is stated only over $\QQ$.  
In that spirit, the 
other assertions are left as an exercise for
the reader.
\end{proof}

To apply Lemma \ref{langex} to $\Pi_1$ and $\Pi_2$, we 
look at (\ref{Pi12}) and label the coefficients 
inside the parentheses as 
$$A_1 = 2 u^p,\,\,\,B_1 = u^{2p} + \frac{du}{c},\,\,\,
A_2 = 2 u^p,\,\,\,
B_2 = u^{2p} + \frac{du}{c} + \frac{1}{4c},$$
so $\Pi_j = T^4 + A_j T^2 + B_j$ modulo
$\kappa^{\times}$-scaling. 
Since we used Lemma \ref{irred4} to prove that each $\Pi_j$ is 
separable and irreducible, we already know that 
$B_j$ and $A_j^2 - 4B_j$ are non-squares in $\kappa(u)^\times$.
A direct calculation also shows that 
$B_j (A_j^2 - 4 B_j)$ is a non-square in $\kappa(u)^\times$.
Therefore, by Lemma \ref{langex}, each of $\Pi_1$ and $\Pi_2$
has Galois group over $\kappa(u)$ that is isomorphic to
$D_4$.
Example \ref{specex}(3) now tells us 
that $\Pi_1$ and $\Pi_2$ 
each have infinitely many irreducible 
$u_0$-specializations.  What about 
simultaneous irreducible specializations?
This is what we need to resolve in order to complete the proof 
of Theorem \ref{phipthm}. 

\begin{theorem}\label{newthm}
There exist infinitely many $u_0$ such that
$\Pi_1(u_0,T)$ and $\Pi_2(u_0,T)$ are both irreducible in
$\kappa_0[T]$.
\end{theorem}

\begin{proof}
Let $L_j/\kappa(u)$ be a splitting field of $\Pi_j$, so
${\rm{Gal}}(L_j/\kappa(u))$ is isomorphic to $D_4$. 
We will show $L_1$ and $L_2$ are linearly disjoint
over $\kappa(u)$.  It will then follow, by 
the Chebotarev density theorem, that 
any pair of Frobenius elements in ${\rm{Gal}}(L_1/\kappa(u)) 
\times {\rm{\Gal}}(L_2/\kappa(u))$ are both attached to infinitely 
many common places on $\kappa(u)$. Theorem \ref{keithfrobold} 
and Example \ref{specex}(3) then imply 
there are infinitely many $u_0$
such that $\Pi_1(u_0,T)$ and $\Pi_2(u_0,T)$ are both irreducible
in $\kappa_0[T]$.

Any intermediate extension in
$L_j/\kappa(u)$, other than $\kappa(u)$, contains a 
quadratic extension of $\kappa(u)$ since 
every proper subgroup of a 2-group is contained in 
a subgroup of index 2. 
We will show that $L_1$ and $L_2$ 
do not contain quadratic subfields (over $\kappa(u)$)
that are $\kappa(u)$-isomorphic, so they must be linearly 
disjoint over $\kappa(u)$.

Inspection shows the only occurrences of non-trivial common factors 
among 
\begin{equation}\label{ablist}
B_1, \,\,\, A_1^2 - 4B_1,\,\,\,B_2,\,\,\,A_2^2 - 4 B_2
\end{equation}
are: the linear polynomial $A_1^2 - 4B_1$ divides $B_1$ and (when
$c = 4 d^{2p}$) the linear polynomial $A_2^2 - 4B_2$ divides $B_1$.
Since $B_1$ is separable with $\deg B_1 > 2$, 
we conclude that the four elements in (\ref{ablist})
are multiplicatively independent modulo squares in
$\kappa(u)^{\times}$. 
This independence modulo squares, coupled with 
the list of quadratic subfields in the $D_4$-case of Lemma \ref{langex}, 
shows $L_1$ and $L_2$ 
do not share a common quadratic extension of $\kappa(u)$.
Thus, they are linearly disjoint 
over $\kappa(u)$.  
\end{proof}

\section{Generic rank bound III. Cohomological arguments}\label{selmer}

By Theorem \ref{phipthm}, there are infinitely many
closed points
$u_0 \in {\mathbf{A}}^1_{\kappa}$
such that the ``specialized'' polynomials
\begin{equation}\label{pi1pi2irred}
\pi_1 = c^{1/p}(T^2 - u_0)^2 + d^{1/p} u_0^{1/p},\,\,\,
\pi_2 = 1 + 4\pi_1 
\end{equation}
in $\kappa_0[T]$ are both irreducible.
{\em These are the only $u_0$ that we shall henceforth consider.}

We view $\pi_1$ and $\pi_2$ as closed points
on ${\mathbf{A}}^1_{\kappa_0} \subseteq \PP^1_{\kappa_0}$. 
The arithmetic of 
\begin{equation}\label{eu00}
\mathscr{E}_{u_0}: y^2 = x^3 + \pi_1^p x^2 - \pi_1^{3p}x
\end{equation}
for such $u_0$ is our focus of interest, as this
will provide the information that we
need to prove that the image of the bottom map
in (\ref{specdig}) has dimension $\le 2$ for these points $u_0$.
This will complete the proof that $\mathscr{E}_{\eta}(F(T))$
has rank 1, thereby concluding the proof of Theorem \ref{mainthm}.

Rather than work with $\mathscr{E}_{u_0}$, it will  
simplify matters to work with the elliptic curve
\begin{equation}\label{pellmodel}
\mathscr{E}'_{u_0}: y^2 = x^3 + \pi_1 x^2 - \pi_1^3 x;
\end{equation}
this elliptic curve is $p$-isogenous to $\mathscr{E}_{u_0} =
(\mathscr{E}'_{u_0})^{(p)}$, so by oddness of $p$
it follows that the map along the bottom of (\ref{specdig})
is canonically identified with the map
\begin{equation}\label{pmap}
\mathscr{E}'_{u_0}(\kappa_0(T))/2 \cdot
\mathscr{E}'_{u_0}(\kappa_0(T)) \rightarrow
\mathscr{E}'_{u_0}(\overline{\kappa}_0(T))/2\cdot
\mathscr{E}'_{u_0}(\overline{\kappa}_0(T)), 
\end{equation}
where $\overline{\kappa}_0$ is an algebraic closure
of $\kappa_0$.   We shall prove that
(\ref{pmap}) is injective for the points $u_0$ presently
under consideration. 

The reduction properties
of the N\'eron model $N(\mathscr{E}_{u_0}) \rightarrow
\PP^1_{\kappa_0}$ were worked out in 
\S\ref{dull} (see above Theorem \ref{phipthm}), and 
the additive and multiplicative properties are
the same for the N\'eron model of
the isogenous elliptic curve $\mathscr{E}'_{u_0}$. 
Thus, letting $j'_{u_0} = j(\mathscr{E}'_{u_0})$, 
for points $u_0$ such that
$\pi_1$ and $\pi_2$ are irreducible in $\kappa_0[T]$ we obtain
the following properties for $N(\mathscr{E}'_{u_0})$: 
\begin{itemize}
\item good reduction at
all closed points of $\PP^1_{\kappa_0}$ away from $\pi_1$ and $\pi_2$,
\item multiplicative reduction
at $\pi_2$, with $\ord_{\pi_2}(j'_{u_0}) = -1$. 
\item additive reduction at $\pi_1$ that is potentially
multiplicative, with $\ord_{\pi_1}(j'_{u_0}) = -2$.
\end{itemize}
By the theory of Tate models for multiplicative reduction,
the component group 
for the N\'eron model at $\pi_2$ is trivial,
so the $\pi_2$-fiber $N(\mathscr{E}'_{u_0})_{\pi_2}$ is a torus.

Fix a geometric point $\overline{\pi}_1$ over the point
$\{\pi_1\} \in \PP^1_{\kappa_0}$.
The reduction at $\pi_1$ is (additive and) potentially
multiplicative, and 
 $\ord_{\pi_1}(j'_{u_0}) = -2$ is negative and even.  We
need to know the structure of
the component group of the additive geometric fiber of the N\'eron
model at $\overline{\pi}_1$.  This can be deduced from Tate's algorithm,
but we give here a direct proof via general principles.

\begin{lemma}\label{compliu}
Let $R$ be a discrete valuation ring with residue field $k$ and fraction
field $K$, and let $E$ be an elliptic curve over $K$ with N\'eron
model $N(E)$ over $R$.  Assume that $\ord_R(j(E))$ is negative and
even, that $E$ has additive reduction over $R$, and that
${\rm{char}}(k) \ne 2$.

If $k$ is perfect and $\overline{k}/k$ is an algebraic closure, then
the geometric component group
$N(E)_{\overline{k}}/N(E)_{\overline{k}}^0$ 
is isomorphic to $\ZZ/2\ZZ \times \ZZ/2\ZZ$.
\end{lemma}

\begin{proof}
The formation of the N\'eron model
over a discrete valuation ring
commutes with base change to a strict henselization 
and to a completion, so we may assume that $k$ is 
separably closed and that $R$ is complete.
Since $\ord_R(j(E)) < 0$, it follows from Tate's theory that
$E$ is a quadratic twist of the Tate curve $E_0$ over $K$
with $j$-invariant $j(E)$.  Since
${\rm{char}}(k) \ne 2$ and
$R$ is strictly henselian, the Tate parameter
$q_{E_0}$ must be a square in $K^{\times}$ because
$\ord_R(q_{E_0}) = -\ord_R(j(E_0)) = -\ord_R(j(E))$ is even.
Thus, the 2-torsion on the
Tate curve $E_0$ is
a constant group over $K$.  This property of the 2-torsion is unaffected
by quadratic twisting, so $E[2]$ is a constant group over $K$.

Using the N\'eron mapping property, we obtain a map of $R$-groups
$$\ZZ/2\ZZ \times \ZZ/2\ZZ \rightarrow N(E),$$
and by passing to $k$-fibers we arrive at a map 
of finite \'etale groups
$$\ZZ/2\ZZ \times \ZZ/2\ZZ \rightarrow N(E)_k/N(E)_k^0.$$
This map is injective by the hypothesis
that the reduction is additive and
${\rm{char}}(k) \ne 2$.
It is a general fact that for any discrete
valuation ring $R$ and any elliptic curve $E$
over the fraction field of $R$, the component group
of the closed fiber of the N\'eron model $N(E)$ 
has order at most 4 when $E$ has additive reduction
and the residue field $k$ is perfect.
This follows from the relationship between
$N(E)$ and the minimal regular proper model
$E^{\rm{reg}}$ over $R$, together with the combinatorial
classification of the extended Dynkin diagrams
that describe the special fiber $E_k^{\rm{reg}}$
(equipped with its intersection form)
when $k$ is {\em algebraically closed};
see \cite[10.2]{liu}.
\end{proof}

\begin{theorem}\label{newcor}
The $2$-torsion subgroup $N(\mathscr{E}'_{u_0})[2]$ is
quasi-finite, \'etale, and separated over $\PP^1_{\kappa_0}$.
It is finite \'etale of order $4$ over $\PP^1_{\kappa_0} - \{\pi_2\}$
and has fiber of order $2$ over $\{\pi_2\}$.
\end{theorem}

\begin{proof}
Since all points on $\PP^1_{\kappa_0}$ have residue
characteristic not equal to 2, doubling
on $N(\mathscr{E}'_{u_0})$ is an \'etale map
that has fiberwise-finite kernel.
Hence, 
$N(\mathscr{E}'_{u_0})[2]$ is a quasi-finite, \'etale, 
and separated $\PP^1_{\kappa_0}$-group, so it
is finite over an open $U \subseteq \PP^1_{\kappa_0}$
if and only if its fiber rank is constant on $U$.
Since $N(\mathscr{E}'_{u_0})_{\pi_2}$ is a torus,
$N(\mathscr{E}'_{u_0})[2]_{\pi_2} = N(\mathscr{E}'_{u_0})_{\pi_2}[2]$ 
has order 2.
For $x \in \PP^1_{\kappa_0} - \{\pi_1, \pi_2\}$
the fiber $N(\mathscr{E}'_{u_0})_x$ is an elliptic curve, 
so its 2-torsion subgroup has order 4. 
It remains to check that 
$N(\mathscr{E}'_{u_0})[2]_{\overline{\pi}_1}$ has order 4.  Consider the
exact sequence of smooth groups
\begin{equation}\label{neronexact}
0 \rightarrow
N(\mathscr{E}'_{u_0})_{\overline{\pi}_1}^0 \rightarrow
N(\mathscr{E}'_{u_0})_{\overline{\pi}_1} \rightarrow
N(\mathscr{E}'_{u_0})_{\overline{\pi}_1}/
N(\mathscr{E}'_{u_0})_{\overline{\pi}_1}^0 \rightarrow 0.
\end{equation}
By Lemma \ref{compliu}, the final term has order 4 and is killed by 2.
Since we are not in characteristic 2, doubling is an automorphism of
the additive group $N(\mathscr{E}'_{u_0})^0_{\overline{\pi}_1}$,
so (\ref{neronexact}) splits.  This gives the result. 
\end{proof}

Consider
the two points 
$P'_0 = (0,0)$
and 
$Q'_0 = (-\pi_1, \pi_1^{2})$
in $\mathscr{E}'_{u_0}(\kappa_0(T))$, 
where $P'_0$ is a rational point of order 2 and 
$(Q'_0)^{(p)} \in (\mathscr{E}'_{u_0})^{(p)}(\kappa_0(T)) =
\mathscr{E}_{u_0}(\kappa_0(T))$ 
is the $u_0$-specialization of (\ref{QdefT}). 

\begin{theorem}\label{comptrick}  
The natural map
\begin{equation}\label{compred}
\mathscr{E}'_{u_0}(\kappa_0(T)) = N(\mathscr{E}'_{u_0})(\PP_{\kappa_0}^1)
\rightarrow
N(\mathscr{E}'_{u_0})_{\overline{\pi}_1}/
N(\mathscr{E}'_{u_0})^0_{\overline{\pi}_1} \simeq 
\ZZ/2\ZZ \times \ZZ/2\ZZ,
\end{equation}
carries $P'_0$ and $Q'_0$ to linearly independent elements.
In particular,
the component group at $\pi_1$ is a constant group 
generated by the classes of $P'_0$ and $Q'_0$.
\end{theorem}

\begin{proof}  
The meaning of the theorem is that
$P'_0$ and $Q'_0$
 reduce into distinct non-identity components
of $N(\mathscr{E}'_{u_0})_{\overline{\pi}_1}$.
By \cite[9.4/35,37]{liu} and \cite[10.2/14]{liu}, 
the smooth locus in
a minimal Weierstrass model is 
the relative
identity component of the N\'eron model
over any discrete valuation ring.  The Weierstrass
model (\ref{pellmodel}) is minimal at $\pi_1$. 
Since $P'_0$ and $Q'_0$ reduce to the unique non-smooth point $(0,0)$ 
on
the closed fiber of this model, 
we conclude 
that the reductions of $P'_0$ and $Q'_0$
in the N\'eron model at $\overline{\pi}_1$ do not lie
in the identity component.

To see that the reductions of $P'_0$ and $Q'_0$ 
in the N\'eron model at $\overline{\pi}_1$ 
lie in distinct components, we just have to check
that the difference $P'_0 - Q'_0 = -(P'_0 + Q'_0)$ also has reduction
not in the identity component on the $\overline{\pi}_1$-fiber of
the N\'eron model; that is,
the point $P'_0 + Q'_0$ should have reduction
$(0,0)$ with respect to the minimal Weierstrass model 
(\ref{pellmodel}) at $\pi_1$.
It is trivial 
to compute
$P'_0 + Q'_0 = (\pi_1^{2}, \pi_1^{3})$,
and this has reduction $(0,0)$.
\end{proof}

\begin{theorem}\label{geominj} Let $\overline{\kappa}_0$ be an algebraic
closure of $\kappa_0 = \kappa(u_0)$,
with $u_0 \in {\mathbf{A}}^1_{\kappa}$ a closed
point such that $\pi_1$ and $\pi_2$ as in $(\ref{pi1pi2irred})$
are irreducible in $\kappa_0[T]$.  
The image of the canonical map
$$c:{\mathscr{E}}'_{u_0}(\kappa_0(T))/2\cdot 
{\mathscr{E}}'_{u_0}(\kappa_0(T))
\rightarrow {\mathscr{E}}'_{u_0}(\overline{\kappa}_0(T))/2\cdot
{\mathscr{E}}'_{u_0}(\overline{\kappa}_0(T))$$
in $(\ref{pmap})$ is spanned by $c(P'_0)$ and $c(Q'_0)$,
so $\dim_{\FF_2} {\rm{image}}(c) \le 2$.
\end{theorem}

\begin{proof}
Let 
$\delta:\mathscr{E}'_{u_0}(\kappa_0(T))/2\mathscr{E}'_{u_0}(\kappa_0(T))
\rightarrow S^{[2]}(\mathscr{E}'_{u_0/\kappa_0(T)})$
be the injective Kummer map to the 2-torsion Selmer group.
Let $K_{x}^{\rm{h}}$ denote the fraction field of
the henselization $\OO_{\PP^1_{\kappa_0},x}^{\rm{h}}$
of the local ring at a closed point $x \in \PP^1_{\kappa_0}$.
For any element ${\sigma} \in S^{[2]}(\mathscr{E}'_{u_0/\kappa_0(T)})$,
the local restriction $${\sigma}_{x} \in {\rm{H}}^1(K_{x}^{\rm{h}},
\mathscr{E}'_{u_0}[2])$$
is in the image of the local Kummer map
$\delta_{x}$ at $x$.  Write ${\sigma}_{x} = \delta_{x}(\xi_x)$ for a point
$$\xi_x \in \mathscr{E}'_{u_0}(K_{x}^{\rm{h}}) =
N(\mathscr{E}'_{u_0})({\mathscr{O}}_{\PP^1_{\kappa_0},x}^{\rm{h}}).$$
By Theorem \ref{comptrick},
subtracting a suitable
$\ZZ$-linear combination of $\delta(P'_0)$ and $\delta(Q'_0)$
from ${\sigma}$ gives a Selmer class ${\sigma}'$ such that
${\sigma}'_{\pi_1} = \delta_{\pi_1}(\xi'_{\pi_1})$, where
$\xi'_{\pi_1}$ reduces into the identity component at $\pi_1$.
Thus,
$S^{[2]}(\mathscr{E}'_{u_0/\kappa_0(T)})$
is generated by $\delta(P'_0)$, $\delta(Q'_0)$, and classes
${\sigma'}$ such that $\sigma'_{\pi_1} = \delta_{\pi_1}(\xi')$
for some local point $\xi'$ in
$\mathscr{E}'_{u_0}(K_{\pi_1}^{\rm{h}})$
that reduces into the identity component
at $\pi_1$; note that this local property of
${\sigma'}$ at $\pi_1$ is independent of the non-canonical choice of
$\xi'$ since any two choices differ
by an element in $[2](\mathscr{E}'_{u_0}(K_{\pi_1}^{\rm{h}}))$
and doubling on $N(\mathscr{E}'_{u_0})_{\pi_1}$
kills the component group (by Lemma \ref{compliu}).

The doubling map on $N(\mathscr{E}'_{u_0})$ is fiberwise surjective
over $\PP^1_{\kappa_0}$ away from $\{\pi_1\}$
and doubling is surjective on the additive identity component
at $\pi_1$ (since $p \ne 2$).  Thus, for
Selmer classes ${\sigma}'$ as above with
local restriction $\sigma'_x = \delta_x(\xi'_x)$,  
the image of 
$\xi'_x$ in $\mathscr{E}'_{u_0}(K_x^{\rm{sh}})$
lies in $[2] {\mathscr{E}}'_{u_0}(K_x^{\rm{sh}})$
for every closed point $x \in \PP^1_{\kappa_0}$ and
every choice of $\xi'_x$
(with $K_{x}^{\rm{sh}}$ denoting a maximal unramified extension
of $K_{x}^{\rm{h}}$).  In other words, the inertial restriction 
$\sigma'_x|_{K_x^{\rm{sh}}}$ is a trivial cohomology
class for all $x$.  Hence,
$S^{[2]}(\mathscr{E}'_{u_0/\kappa_0(T)})$
is spanned by the images of $P'_0$ and $Q'_0$ and
the intersection of this Selmer group with
the subgroup of everywhere unramified classes in
${\rm{H}}^1(\kappa_0(T),\mathscr{E}'_{u_0}[2])$.

Let us
now recall how to describe the group of everywhere
unramified classes in terms of \'etale cohomology.
Let $G = N(\mathscr{E}'_{u_0})[2]$ and $\mathbf{P} = \PP^1_{\kappa_0}$, so
$G$ is a quasi-finite separated \'etale commutative ${\mathbf{P}}$-group.
If we let $i_{\eta}:\eta \rightarrow {\mathbf{P}}$ be the canonical
map from the generic point $\eta$ of ${\mathbf{P}}$, then the identity
$N(\mathscr{E}'_{u_0}) = i_{\eta \ast}(\mathscr{E}'_{u_0})$
on the smooth site over $\PP$ 
implies $G = i_{\eta \ast}(G_{\eta})$ as
\'etale sheaves (by passing to 2-torsion subsheaves).  Thus,
using the \'etale topology, 
the Leray spectral sequence
${\rm{E}}_2^{r,s} = 
{\rm{H}}^r({\mathbf{P}}, {\rm{R}}^s i_{\eta \ast}(G_{\eta}))
\Rightarrow {\rm{H}}^{r+s}(\eta,G_{\eta})$
has ${\rm{E}}_2^{r,0} = {\rm{H}}^r({\mathbf{P}},G)$, so we get
an exact sequence of low-degree terms
\begin{equation}\label{leray}
0 \rightarrow {\rm{H}}^1({\mathbf{P}},G) 
\stackrel{\alpha}{\rightarrow}
{\rm{H}}^1(\eta,G_{\eta}) \stackrel{\oplus \beta_x}{\rightarrow}
\bigoplus_x {\rm{H}}^0(\kappa_0(x),{\rm{H}}^1(K_x^{\rm{sh}},G)).
\end{equation}
Here $\alpha$ is the canonical restriction map
to the generic point and $\beta_x$ is the canonical local
restriction map at the non-generic point
$x$ of ${\mathbf{P}}$.
Hence, ${\rm{H}}^1({\mathbf{P}},G) \subseteq {\rm{H}}^1(\eta,G_{\eta})$
is the group of everywhere unramified classes.

In view of the preceding considerations, to prove Theorem 
\ref{geominj} it suffices to prove that the restriction map
$${\rm{H}}^1(\kappa_0(T),\mathscr{E}'_{u_0}[2]) \rightarrow
{\rm{H}}^1(\overline{\kappa}_0(T),\mathscr{E}'_{u_0}[2])$$
kills the subgroup ${\rm{H}}^1({\mathbf{P}},G)$
of everywhere unramified classes, where
$G = N(\mathscr{E}'_{u_0})[2]$.  
We will 
prove the stronger assertion that the map 
${\rm{H}}^1(\PP,G) \rightarrow {\rm{H}}^1(\PP_{\overline{\kappa}_0},G)$
vanishes.

Let $U' = {\mathbf{P}} - \{\pi_2\}$ and let $j':U' \hookrightarrow 
{\mathbf{P}}$
be the canonical open immersion.  By Theorem \ref{newcor},
$G|_{U'}$ is finite \'etale over $U'$ and 
$G_{\pi_2}$ has order 2 over $\kappa(\pi_2)$.  Thus, 
the nontrivial 2-torsion point $(0,0)$ defines a short
exact sequence of \'etale sheaves
\begin{equation}\label{filterg}
0 \rightarrow \ZZ/2\ZZ \rightarrow G \rightarrow j'_{!}(\ZZ/2\ZZ)
\rightarrow 0
\end{equation}
over ${\mathbf{P}}$.  By considering the exact sequence
of pullback sheaves on $\PP_{\overline{\kappa}_0} =
\PP^1_{\overline{\kappa}_0}$ and using the vanishing
of ${\rm{H}}^1(\PP^1_{\overline{\kappa}_0},\ZZ/2\ZZ)$, we arrive at
a commutative square
\begin{equation}\label{h1dig}
\xymatrix{
{{\rm{H}}^1(\PP_{\overline{\kappa}_0},G)} \ar[r] &
{{\rm{H}}^1(\PP_{\overline{\kappa}_0},j'_{!}(\ZZ/2\ZZ))}\\
{{\rm{H}}^1(\PP,G)} \ar[u] \ar[r] & {{\rm{H}}^1(\PP,j'_{!}(\ZZ/2\ZZ))}
\ar[u] }
\end{equation}
whose top side is injective
and whose vertical maps are 
the natural pullback maps. 
It therefore suffices to prove that the pullback map on the right
side vanishes, 
and this is a property that does not involve $G$.

Since $j'_{!}(\ZZ/2\ZZ)$ is
represented by an \'etale $\PP$-group that is quasi-affine over
$\PP$ (it is the complement in $(\ZZ/2\ZZ)_{\PP}$
of the non-identity point over $\{\pi_2\} \in \PP$), 
the elements of ${\rm{H}}^1(\PP,j'_{!}(\ZZ/2\ZZ))$
are in bijection with isomorphism classes of (representable)
\'etale $j'_{!}(\ZZ/2\ZZ)$-torsors
on $\PP$; the same holds over $\PP_{\overline{\kappa}_0}$,
 and the right side of (\ref{h1dig})
is thereby identified with base-change on
torsors.  Thus, we just need
to prove that every \'etale $j'_{!}(\ZZ/2\ZZ)$-torsor
on $\PP$ has a $\PP_{\overline{\kappa}_0}$-point.

Let $i:\Spec \kappa_0(\pi_2) \hookrightarrow \PP$
be the closed complement to $U'$, so we have a short exact sequence
\begin{equation}\label{filterz2}
0 \rightarrow j'_{!}(\ZZ/2\ZZ) \rightarrow \ZZ/2\ZZ \rightarrow
i_{\ast}(\ZZ/2\ZZ) \rightarrow 0
\end{equation}
of \'etale sheaves on $\PP$, with $i_{\ast}(\ZZ/2\ZZ)$ supported
at {\em one} physical point on $\PP$.  Thus,
the natural map 
$${\rm{H}}^1({\mathbf{P}},j'_{!}(\ZZ/2\ZZ)) \rightarrow
{\rm{H}}^1({\mathbf{P}},\ZZ/2\ZZ)$$
is injective.
Since  
${\rm{H}}^1({\mathbf{P}},\ZZ/2\ZZ) = {\rm{H}}^1(\PP^1_{\kappa_0},\ZZ/2\ZZ)
= {\rm{H}}^1(\kappa_0,\ZZ/2\ZZ)$, 
clearly ${\rm{H}}^1(\PP,\ZZ/2\ZZ)$
has order 2 with its nontrivial element represented by
the nontrivial $\ZZ/2\ZZ$-torsor $\PP_{\kappa'_0} \rightarrow \PP$
for a quadratic extension $\kappa'_0/\kappa_0$,
and the subgroup ${\rm{H}}^1(\PP,j'_{!}(\ZZ/2\ZZ))$
has order 1 or 2. 

The fiber of the $\ZZ/2\ZZ$-torsor $\PP_{\kappa'_0}$
over the point $\pi_2 \in \PP$ 
is a {\em split} double cover of
$\Spec \kappa_0(\pi_2)$ since $\deg_{\kappa_0}\pi_2 = 4$ and $[\kappa'_0:\kappa_0] = 2$. 
Removing one of the two points over $\{\pi_2\}$
in $\PP_{\kappa'_0}$ gives an open
${\mathscr{T}} \subseteq \PP_{\kappa'_0}$ that is a 
nontrivial $j'_{!}(\ZZ/2\ZZ)$-torsor over $\PP$.
Hence, ${\rm{H}}^1(\PP,j'_{!}(\ZZ/2\ZZ))$ has order 2 and 
contains ${\mathscr{T}}$ as its unique nontrivial element.
Since ${\mathscr{T}}$ obviously acquires a section upon extending
the ground field $\kappa_0$ to $\kappa'_0$, we
conclude that ${\mathscr{T}}(\PP_{\overline{\kappa}_0})$
is nonempty. 
\end{proof}

\section{Nagao's function field conjecture}\label{nagsec}

We have already shown  
that, under the parity conjecture, (\ref{mainfamily}) 
has ``unexpected'' rank behavior for its specializations in $T$.
A conjecture of Nagao~\cite[p.~14]{nag2}, to be reviewed below, predicts 
that the Mordell--Weil 
rank of a non-constant elliptic curve over $\QQ(T)$ 
is a certain limit of averages over mod $p$ specializations 
as $p$ runs over the primes.  Nagao's conjecture
admits a natural variant for elliptic curves over
$F(T)$ for any global field $F$, and so in view of the
unusual behavior
of Mordell--Weil ranks in (\ref{mainfamily}) discussed
in this paper
one may be led to wonder if it is necessary to introduce
restrictions in the analogue of Nagao's
conjecture over $F(T)$ for global function fields $F$.
In this final section we shall check
numerically that (\ref{mainfamily}) appears to satisfy
a straightforward $\kappa(u)(T)$-analogue of Nagao's conjecture,
and so the surprising behavior of fibral Mordell--Weil ranks
in (\ref{mainfamily}) (conditional on the parity conjecture) does
not seem to cast doubt on the standard reformulation of
Nagao's conjecture in the case of pencils
of elliptic curves over global function fields.

Nagao's conjecture over $\QQ(T)$ involves the following data. 
Pick a non-constant elliptic curve ${\mathscr E}_\eta$
over $\QQ(T)$, a prime $p$, 
and let ${\mathscr E}_p$ be the reduction modulo $p$. 
For all but finitely many $p$, ${\mathscr E}_p$ is 
an elliptic curve over $\FF_p(T)$.
Furthermore,
the fiber ${\mathscr E}_{p,s}$ 
over each $s \in \PP^1(\FF_p)$ is an elliptic curve 
over $\FF_p$ except for a set of $s$ with size bounded 
uniformly in $p$. Set
\begin{equation}\label{Apdefn}
A_p({\mathscr E}_\eta) = 
\frac{1}{p}\sum_{s \in \PP^1(\FF_p)}
a_p({\mathscr E}_{p,s}), 
\end{equation}
where we only include points $s$ for which
the reduction $\mathscr{E}_{p,s}$ over $\FF_p$
is an elliptic curve, with $a_p$ denoting the Frobenius trace for
such terms.  
Nagao conjectured that
the rank of ${\mathscr E}_\eta(\QQ(T))$ can be computed as 
\begin{equation}\label{nagcon}
{\rm rank}({\mathscr E}_\eta(\QQ(T)))
\stackrel{?}{=} 
\lim_{n \rightarrow \infty} 
\frac{1}{n}\sum_{p \leq n}-A_p({\mathscr E}_\eta)\log p.
\end{equation}
The existence of the limit is part of the conjecture, and 
bad $p$ are omitted from the sum.
(In the definition (\ref{Apdefn}), the omission of a uniformly bounded
number of $s$ for each $p$ is irrelevant for (\ref{nagcon}).
Moreover, division by $p$ rather than by $p+1 = \#\mathbf{P}^1(\FF_p)$
also has no effect because
$|a_p| \leq 2\sqrt{p}$, by the Riemann hypothesis, and 
$(1/n) \sum_{p \le n} (\log p)/\sqrt{p} \sim 2/\sqrt{n} = o(1)$.)
Rosen and Silverman~\cite{rossil} have shown the truth of 
(\ref{nagcon}) is 
closely related to a conjecture of Tate.
Nagao actually made two conjectures about 
formulas for ranks, one over $\QQ$~\cite[p.~213]{nag1} and the other 
over $\QQ(T)$ as above.  The latter was inspired by the former.

The analogue of Nagao's conjecture for an elliptic curve 
${\mathscr E}_\eta$ 
over $\kappa(u)(T)$ ought to be the following.  
For a place $v$ on $\kappa(u)$, with residue field $\FF_v$, set 
\begin{equation}\label{Avdefn}
A_v({\mathscr E}_\eta) = 
\frac{1}{\Nm\!v}\sum_{s \in \PP^1(\FF_v)}
a_v({\mathscr E}_{v,s})
\end{equation}
when $\mathscr{E}_{\eta}$ has good reduction
over $\FF_v(T)$, and we only include
$s$ such that the reduction
$\mathscr{E}_{v,s}$ at $s$ is an elliptic curve.  
Here $\Nm\!v = \#\FF_v$ and 
$a_v$ is the Frobenius trace.  Then, setting $q = \#\kappa$,  
the conjecture is
\begin{eqnarray*}
{\rm rank}({\mathscr E}(\kappa(u)(T)))
& \stackrel{?}{=} &
\lim_{n \rightarrow \infty} 
\frac{1}{q^n}\sum_{\deg v = n}-A_v({\mathscr E}_\eta)\deg v \\
& = & 
\lim_{n \rightarrow \infty} 
\frac{n}{q^n}\sum_{\deg v = n}-A_v({\mathscr E}_\eta).
\end{eqnarray*}
(Allowing all $v$ with $(\deg v)|n$ in the first
sum adds a term $O(n/q^{n/2}) = o(1)$, so whether we sum
over those $v$ satisfying $\deg v = n$ or $(\deg v)|n$ is irrelevant
for the meaning of this conjecture over $\kappa(u)$.)
As in the situation over $\QQ$,
the Riemann hypothesis for elliptic curves over finite fields ensures that
this conjecture is unaffected by omitting terms 
in the $A_v({\mathscr E}_{\eta})$'s
corresponding to a set of 
points $s \in \PP^1(\FF_v)$ whose size is bounded
uniformly in $\deg v$.

We now consider this conjecture numerically for 
the curve in (\ref{mainfamily}) where 
$c = d = 1$. (Then $h(T) = T^{2p} + u$ in
(\ref{hdef}).)
The generic rank is 1 by Theorem \ref{mainthm}, so we expect
\begin{equation}\label{qqq}
\sum_{\deg v = n} -A_v({\mathscr E}_\eta) \stackrel{?}{\sim} \frac{q^n}{n}.
\end{equation}
From the definition of $A_v$ in (\ref{Avdefn}), 
clearly (\ref{qqq}) is the same as
\begin{equation}\label{asymp}
\sum_{\deg v = n}\sum_{s \in \PP^1(\FF_v)}
-a_v({\mathscr E}_{v,s}) \stackrel{?}{\sim} \frac{q^{2n}}{{n}}.
\end{equation}
Consider only $n \geq 2$, since in degree 1 
the specialization 
${\mathscr E}_{\infty}$ from (\ref{mainfamily}) 
is not an elliptic curve over $\kappa(T)$.
Writing $a_v$ in terms of the quadratic character 
on $\FF_v$, an equivalent reformulation of (\ref{asymp}) is
\begin{equation}\label{last}
\sum_{\deg \pi = n}
\sum_{s \in \kappa[u]/(\pi)}
\sum_{x \in \kappa[u]/(\pi)}
\left(\frac{x^3+h(s^2+u)x^2-h(s^2+u)^3x}{\pi}\right) 
\stackrel{?}{\sim} \frac{q^{2n}}{n}, 
\end{equation}
where $\pi$ runs over monic irreducibles in $\kappa[u]$ of degree $n$ 
and $(\frac{\cdot}{\pi})$ is a Legendre symbol.
(Technically, we should be omitting from (\ref{last}) the at most 8 values of 
$s$ in each $\FF_v$ such that ${\mathscr E}_{v,s}$ is not an elliptic 
curve, but the inclusion of such terms in (\ref{last}) contributes 
an amount which has strictly smaller growth than $q^{2n}/n$, so there 
is nothing lost by allowing such $s$ for the sake of a cleaner summation.)

Tables \ref{tab3data} and \ref{tab5data} 
compare the two sides of (\ref{last}) 
for $q = 3$ and $q = 5$. 
All decimal approximations are truncated after the third 
digit beyond the decimal point. 
This evidence supports the truth
of Nagao's conjecture for (\ref{mainfamily}).  Thus, 
we do not find any evidence that the formulation
of Nagao's conjecture for pencils over global
function fields requires unexpected restrictions. 
(We thank Aaron Silberstein for his assistance with the 
preceding numerical calculations.)

\begin{table}
\begin{center}
\begin{tabular}{c|r|r|r} 
$n$ & Left side &  Right side & Ratio \\ \hline
2 & 17 & 40.500& 2.382\\
3 & 173 & 243.000& 1.404\\
4 & 1186 & 1640.250& 1.383\\
5 & 10788& 11809.800& 1.094\\
6 & 91816& 88573.500& .964\\
\end{tabular}
\caption{(\ref{last}) with $q = 3$}\label{tab3data}
\end{center}
\end{table}

\begin{table}
\begin{center}
\begin{tabular}{c|r|r|r} 
$n$ & Left side &  Right side & Ratio \\ \hline
2 &228 &312.500 & 1.370\\
3 &5430 &5208.333 & .959\\
4 &96802 &97656.250 & 1.008\\
\end{tabular}
\caption{(\ref{last}) with $q = 5$}\label{tab5data}
\end{center}
\end{table}

\appendix

\section{Known results over $\QQ$}\label{avgrtnum}

In the Introduction, we saw how to search for (conditional)
examples of elevated rank over $\QQ(T)$: assume the 
parity conjecture over $\QQ$ and 
try to construct an elliptic curve over
$\QQ(T)$ that satisfies (\ref{Wnot}) for all but finitely many $t \in 
\PP^1(\QQ)$. 
We wish to explain why this sufficient strategy
is essentially
necessary if
we also assume three additional standard conjectures
over $\QQ$.
Moreover, we will see that if all of these conjectures are true then
there do not exist
non-isotrivial examples of elevated rank over $\QQ(T)$. 
The three additional conjectures we bring in are: the density 
conjecture, the squarefree-value conjecture, and Chowla's 
conjectures.

The {\it density conjecture} says that 
for any elliptic curve $\mathscr{E}_{\eta}$ over $\QQ(T)$,
the rank of $\mathscr{E}_t(\QQ)$ 
equals 
${\rm{rank}}(\mathscr{E}_{\eta}(\QQ(T)))$ or 
${\rm{rank}}(\mathscr{E}_{\eta}(\QQ(T))) + 1$ 
except for a set of $t \in \PP^1(\QQ)$ with density 0, as measured by 
height.  Granting this and the parity conjecture, any example of 
elevated rank over $\QQ(T)$ will satisfy (\ref{Wnot}) 
for all $t$ outside of a set of density 0. 
Therefore, if $\mathscr{E}_{\eta}$ has elevated
rank then the average value of $W({\mathscr E}_t)$,
in the sense of the following definition, 
is either 1 or $-1$.

\begin{definition}
For any elliptic curve $\mathscr{E}_{\eta}$ 
over $\QQ(T)$, its {\em average root number} is
\begin{equation}\label{eq:nothung}
\avg_{\QQ} W(\mathscr{E}_t)
:= \lim_{N\to \infty}
\frac{\sum_{t\in \PP^1(\QQ), h_{\QQ}(t)<N} W(\mathscr{E}_t)}{\#\{t\in
\PP^1(\QQ), h_{\QQ}(t)<N\}}
\end{equation}
if this limit exists, 
where $h_{\QQ}$ is the standard logarithmic height function on 
$\PP^1({\QQ})$ (defined by the standard normalized collection of
absolute values on $\QQ$).
In the summation,   
the finitely many $t$ at which $\mathscr{E}_t$ is non-smooth
are dropped out. 
\end{definition}

The existence
of the average root number is not evident {\em a priori}, and 
its value might depend on the choice of coordinate on $\PP^1$.
(The height $h_{\QQ}$ depends on the coordinate.)
If the average exists, then clearly 
$-1 \leq \avg_\QQ W(\mathscr{E}_t) \leq 1$.  
If we assume the parity and density conjectures, then 
any example of elevated rank over $\QQ(T)$ must have 
average root number 1 or $-1$. 

\begin{remark} 
For any elliptic curve $E_0$ over $\QQ$, 
Rizzo~\cite{rizzo1} proved that the set of average
root numbers that unconditionally exist
for quadratic
twists of $E_0$ over $\QQ(T)$
is dense in the interval $[-1,1]$. 
\end{remark}

We now introduce the squarefree-value conjecture and Chowla's conjectures;
these lead to a formula for 
$\avg_\QQ W(\mathscr{E}_t)$ for any 
elliptic curve $\mathscr{E}_{\eta}$ over $\QQ(T)$. 
This formula turns out never to equal 1 or $-1$
for non-isotrivial elliptic curves 
over $\QQ(T)$, thereby (conditionally) ruling out
the possibility of elevated rank for such elliptic curves.

The {\em squarefree-value conjecture}
says that a polynomial over $\ZZ$ takes squarefree values as often as
is suggested by naive 
probabilistic heuristics.  For example, 
if $f(T) \in \ZZ[T]$ is squarefree, then the prediction is
$$
\#\{1 \leq n \leq x : f(n) \text{ is squarefree}\} 
\sim Cx 
$$
as $x \rightarrow \infty$, 
where $C = \prod_p \left(1 - c_p/p^2\right)$
with $c_p$ denoting the number of solutions to $f(T) = 0$ in $\ZZ/(p^2)$. 
(If $1 - c_p/p^2 = 0$ for some $p$, then $C = 0$ and obviously 
$f(n)$ is never squarefree.  Otherwise $C$ is an
absolutely convergent (positive) product.) 
We refer the reader to work of Granville~\cite{granville} 
for a more complete statement of this conjecture, 
including the variant for homogeneous polynomials in two variables over
$\ZZ$.  The squarefree-value conjecture is known unconditionally 
for polynomials in $\ZZ[T]$ 
with small degree, and Granville~\cite{granville} deduced 
the general case (all degrees) from the {\it abc}-conjecture.  
(Poonen~\cite{poonen} extended these results to polynomials in 
any number of variables over $\ZZ$, but only the cases 
treated by Granville in one and two variables are related to the 
variation of root numbers in pencils of elliptic curves over $\QQ$.)

\begin{remark} Low-degree proved instances of 
the squarefree-value conjecture were used in the 
study of ranks of elliptic curves over $\QQ$ in 
\cite{gouveamazur}, where 
families of quadratic twists were considered. 
\end{remark}

The final conjecture we need over $\QQ$, 
due to Chowla \cite[p.~96]{chowla}
in the one-variable case, 
concerns the average behavior 
of the Liouville 
function on values of a polynomial. 
Recall that Liouville's function $\lambda$ 
is the totally multiplicative function 
on $\ZZ$ defined by $\lambda({\pm}p) = -1$ when $p$ is prime, 
$\lambda(\pm 1) = 1$, and $\lambda(0) = 0$. 

The {\em one-variable Chowla conjecture} says 
that for any non-constant $f(T)$ in $\ZZ[T]$, 
which is not a perfect square up to sign, 
the sequence $\lambda(f(n))$ has average value 0 as $n$ runs over 
any arithmetic progression.  That is, for any 
arithmetic progression $a+b\ZZ$ (with $a \in \ZZ$ and $b \in \ZZ$, 
$b \ne 0$), as $N \rightarrow \infty$ we have
$$
\frac{\sum_{n \in (a + b\ZZ) \cap [0,N]} \lambda(f(n))}{
\#((a + b\ZZ) \cap [0,N])} \rightarrow 0.
$$
(Clearly, 
with a linear change of variables, we can 
state this as a conjecture over all $f$ using 
only $a = 0$ and $b = 1$.  We prefer the above 
superficially more general form because it 
matches the two-variable 
conjecture more closely.)

The {\em two-variable Chowla conjecture} says that 
for any non-constant homogeneous $f$ in $\ZZ[U,V]$ which is not a 
perfect square up to sign,  
the sequence $\lambda(f(m,n))$ has average value 0 as $(m,n)$ runs over
lattice points in any sector of the plane with vertex
at the origin.  More precisely, for any
coset $L \subseteq \ZZ^2$ of 
an arbitrary sublattice of $\ZZ^2$, and any
open sector 
$S \subseteq \RR^2$ with positive 
angular measure
and vertex at the origin,
\begin{equation}\label{eq:chow2}
\frac{\sum_{(m,n) \in S \cap L \cap [-N,N]^2} \lambda(f(m,n))}{
\#(S \cap L \cap [-N,N]^2)} \rightarrow 0
\end{equation}
as $N \rightarrow \infty$. 
If the condition $(m,n) = 1$ is imposed on the terms in the sum in 
(\ref{eq:chow2}),
then the resulting general conjecture 
is logically equivalent to
the general conjecture (\ref{eq:chow2}).

In \cite{helfgott} and \cite{helfgott2}, 
the squarefree-value conjecture and the two-variable 
Chowla conjecture are used to derive 
(conditional) formulas for $\avg_{\QQ} W(\mathscr{E}_t)$
for any $\mathscr{E}_{\eta/\QQ(T)}$.  
The analysis falls into two cases: 

\begin{itemize}

\item
[Case 1:] 
\label{first}
The minimal regular proper model 
$\mathscr{E} \rightarrow \PP^1_{\QQ}$ 
has no nodal geometric fiber.  (That is, $\mathscr{E}_{\eta}$
has no 
point of multiplicative reduction on $\PP^1_{\QQ}$.) 

\item 
[Case 2:]
\label{second}
The fibration $\mathscr{E} \rightarrow \PP^1_{\QQ}$
has a nodal geometric fiber.

\end{itemize}

Consider a non-isotrivial
$\mathscr{E}_{\eta/\QQ(T)}$ in Case 1.  
Let $M_t$ denote the finite set of primes $p \in \ZZ$ such that
${\mathscr E}_t$ has multiplicative reduction
at $p$.
The collection $\{M_t\}_{t \in \PP^1(\QQ)}$ is restricted
as $t$ 
varies, in the following sense.
Assuming the square-free
value conjecture, we have that, for any small $\varepsilon > 0$, 
there is a finite set of 
prime numbers 
$\mathcal{S}_{\varepsilon}$ such that the set of 
$t \in \PP^1(\QQ)$ with $M_t \subseteq \mathcal{S}_{\varepsilon}$
has height density $\ge 1 - \varepsilon$.
That is, 
roughly speaking, ``most'' fibers have 
their primes of multiplicative reduction lying in a 
common finite set.
(This remark is implicit in
\cite[Lemma~2.1]{man}.)
Moreover, for such $t$ the bad primes for
$\mathscr{E}_{t/\QQ}$ outside of $\mathcal{S}_{\varepsilon}$ are 
the prime factors of values of
certain irreducible primitive polynomials over $\ZZ$
that correspond to the points of additive reduction for 
$\mathscr{E}_{\eta}$ on $\PP^1_{\QQ}$. 
(In 
particular, these primitive polynomials are independent of $t$
and $\varepsilon$.)
For the study of average
root numbers of  elliptic curves 
over $\QQ$, the essential difference between additive and multiplicative 
reduction is the simpler statistical variation 
for local root numbers in the additive case.  (See the formulas in Theorem 
\ref{localrtnum}.) 
Assuming the squarefree-value conjecture,
for ``most'' $t$ 
the set of bad primes for $\mathscr{E}_t$ outside of 
$\mathcal{S}_{\varepsilon}$ 
can be controlled, and a formula 
\begin{equation}\label{avgform}
\avg_{\QQ} W(\mathscr{E}_t) = C_{\infty}
\prod_{p} C_p
\end{equation} 
is thereby obtained, 
where $C_{\infty}$ is an algebraic number in $\RR$, 
each $C_p$ is a non-zero rational number, and $\prod_p C_p$ is
an absolutely convergent (non-zero) product.

Here are two examples that 
illustrate (\ref{avgform}) (not elevated rank) for non-isotrivial 
elliptic curves in Case 1. 

\begin{example}\label{washex}
An example of Washington \cite[\S3]{washington}
over $\QQ(T)$ is 
\begin{equation}\label{washcurve}
\mathscr{E}_{\eta}:y^2 = x^3 + Tx^2 - (T+3)x + 1.
\end{equation}
The point $(0,1)$ on this curve
has infinite order (use 
Theorem \ref{nlutzFT}, as in the proof
of Corollary \ref{rankup}).  Since
$\mathscr{E}$ is a rational surface, 
it is not difficult to prove (using either analytic methods
of Rosen--Silverman or a reduction to positive characteristic 
and algebraic methods 
of Artin--Tate) that
$\mathscr{E}_{\eta}(K(T))$ has rank 1 for every number field $K$.

In \cite{rizzo2}, Rizzo shows
$W({\mathscr{E}}_t) = -1$ for every $t \in \ZZ$. 
However, $W(\mathscr{E}_t) = 1$ for many non-integral $t \in \QQ$, 
such as (using PARI) $t = -1/2, 1/3$, and $3/2$.
An application of one of the proved instances of 
the squarefree-value conjecture in low degree
shows that (\ref{avgform}) is unconditionally true 
for $\mathscr{E}_{\eta}$ in (\ref{washcurve}), and a computation yields 
$C_{\infty} = 0$.  Therefore, in this example, 
 $\avg_{\QQ} W(\mathscr{E}_t) = 0$ unconditionally. 
\end{example}

\begin{example}
Let $f(T) = -5 - 2 T^2$ and $g(T) = 2 + 5 T^2$. Consider
$${\mathscr{E}}_{\eta}:y^2 = x^3 + a(T)x + b(T)$$
over $\QQ(T)$, where 
$$a(T) = -27f g(f^3 - g^3)^2,\,\,\,
b(T) = -\frac{54(f^3 + g^3)(f^3 - g^3)^3}{2}.$$
Low-degree proved
instances of the squarefree-value conjecture imply
that the conditional formula (\ref{avgform}) is 
true for this $\mathscr{E}_{\eta}$.  This leads to the 
explicit formula
$$
\avg_{\mathbf{Q}} W({\mathscr{E}}_t) = 
\frac{1}{6}\cdot \prod_{p \ne 2, 3, 7, 19}
\left(1 - \frac{a_p}{(p+1)^2}\right) = 0.1527\dots,
$$
where $a_p = 1 + \chi_p(-1) + (1 + \chi_3(p))(1 + \chi_{19}(p))$
and $\chi_{\ell}$ is the mod $\ell$ Legendre symbol.
\end{example}

A closer analysis of the work that leads to (\ref{avgform}) 
in Case 1 shows 
that if the squarefree-value
conjecture is assumed then $\avg_{\QQ} W(\mathscr{E}_t)$ 
{\em cannot} equal 
$1$ or $-1$ in Case 1 when ${\mathscr{E}}_{\eta}$ is non-isotrivial. 
Therefore, if the density conjecture, parity conjecture, 
and squarefree-value conjecture are true, then
in Case 1 there does not exist a non-isotrivial 
elliptic curve over $\QQ(T)$ with elevated rank.

We turn now to Case 2, 
so $\mathscr{E} \rightarrow \PP^1_{\QQ}$ 
has a nodal geometric fiber.  
Such $\mathscr{E}$ must be non-isotrivial.
The reasoning in Case 1 breaks down, since 
there do not exist sets of $t$ with height density 
arbitrarily close to 1 such that 
the ${\mathscr{E}_t}_{/\QQ}$'s have multiplicative reduction in 
a common finite set of primes. 
Now there is a  non-constant homogeneous two-variable
polynomial $f_{\mathscr{E}} \in \ZZ[U,V]$, which is not a square, 
such that as $t \in \PP^1(\QQ)$ varies with $\mathscr{E}_{t/\QQ}$ smooth,
the 
variation of the product of the local root numbers of $\mathscr{E}_{t/\QQ}$
at the places of multiplicative reduction 
is governed by the variation of 
$\lambda(f_{\mathscr{E}}(m,n))$
where $m/n$ is the reduced form of $t$.
A similar phenomenon 
happens in our function field examples, 
using the Liouville function on $\kappa[u]$ that
assigns value $-1$ to irreducibles and 
extends to all of $\kappa[u]$ 
by total multiplicativity; see 
Remark \ref{remlambda}.) 
If we assume Chowla's two-variable conjecture for
$f_{\mathscr{E}}$, then 
the variation of $\lambda(f_{\mathscr{E}}(m,n))$ 
as $t = m/n$ varies can be controlled. 
Using this, in \cite[\S1.7]{helfgott} it is shown that 
if the squarefree-value conjecture
is also assumed, then 
$\avg_{\QQ} W(\mathscr{E}_t)$ exists and equals 0;
in particular, this average does not equal 1 or $-1$. 
Thus, the parity, density, squarefree-value, and Chowla 
conjectures predict 
that no elliptic curve in Case 2 
has elevated rank.

Our discussion of $\avg_{\QQ} W(\mathscr{E}_t)$ has shown that if we 
assume the squarefree-value conjecture and Chowla's 
two-variable conjecture then 
this average cannot equal 1 or $-1$ if
$\mathscr{E}_\eta$ is non-isotrivial.
If we accept the parity conjecture 
and the density conjecture, 
then any example of elevated rank over $\QQ(T)$ 
has $\avg_{\QQ} W({\mathscr E}_t) = 1$ or $-1$. 
Therefore, if all four conjectures are true then 
all examples of elevated rank over 
$\QQ(T)$ must be isotrivial.

\section{The surprise in characteristic $p$}\label{contrsec}

We now replace $\QQ$ with $F = \kappa(u)$
and replace $\ZZ$ with $\kappa[u]$, where $\kappa$ is any finite field. 
For the moment, $\kappa$ may have characteristic 2.
Granting the parity conjecture over $F$, 
no new ideas should be required
to construct isotrivial examples of elevated rank over $F(T)$ 
analogous to the examples of Cassels--Schinzel and Rohrlich. 
(The case of characteristic 2 is presumably more delicate.) 
We want to explain why it is reasonable
to expect {\em a priori} that non-isotrivial examples
of elevated rank might exist over $F(T)$, despite
the conclusions over $\QQ(T)$ in Appendix \ref{avgrtnum}.

The squarefree-value conjecture for multivariable polynomials
over $\kappa[u]$, for $\kappa$ with any characteristic, 
was proved by Ramsay \cite{ramsay} in the separable case
for one variable, and was proved by Poonen \cite{poonen} in general.
Thus,
provided that ${\rm{char}}(F) \ne 2, 3$ (to avoid
problems with wild ramification at arbitrarily many places), 
the methods used over $\QQ(T)$ can be adapted 
to prove
an {\em unconditional} formula akin to (\ref{avgform})
in the analogue of 
Case 1 in Appendix \ref{avgrtnum}.
(This 
is the case of elliptic curves $\mathscr{E}_{\eta/F(T)}$  
such that $\mathscr{E} \rightarrow \PP^1_F$ does not have
any nodal geometric fibers.)  However, 
to adapt the $\QQ(T)$-methods
to prove that $\avg_F W(\mathscr{E}_t)$ is strictly between
$1$ and $-1$ for a non-isotrivial
$\mathscr{E} \rightarrow \PP^1_F$ without nodal fibers,
we need to impose a restriction that is
always satisfied in characteristic 0: the points in the
(non-empty) support of the conductor of $\mathscr{E}_{\eta}$
on $\PP^1_F$ are \'etale over $F$.  We expect
that if this \'etale restriction
on the support of the conductor is dropped,
then there should be non-isotrivial examples without
nodal geometric fibers such that 
the average root number is $1$ or $-1$.
Moreover, in all positive characteristics there should exist such examples
that also have elevated rank (granting the parity conjecture).

Let us now turn to the analogue of Case 2 from Appendix
\ref{avgrtnum}, so $\mathscr{E} \rightarrow \PP^1_F$
has some nodal geometric fibers.  The study of such 
elliptic fibrations 
in characteristic 0 uses
Chowla's conjectures over $\ZZ$,
as we saw in Appendix \ref{avgrtnum}.
However, there are counterexamples to the $\kappa[u]$-analogues of 
Chowla's conjectures. 
In \cite{ccg}, it is shown that counterexamples
to Chowla's one-variable conjecture
are a common (but not ``generic'')
phenomenon.  
For example, elementary (but non-obvious) methods show that
for any finite field $\kappa$ with arbitrary characteristic $p$, 
$f(T) = T^{4p} + u \in \kappa[u][T]$
violates the one-variable Chowla conjecture: $\lambda(f(g)) = 1$ 
for every $g \in \kappa[u]$ with $g \not\in \kappa$.
Similarly, in the sense of Chowla's
two-variable conjecture, the homogeneous polynomial 
$$
aX^{4p} + buY^{4p} \in \kappa[u][X,Y]
$$
with $a, b \in \kappa^{\times}$ has rather non-random
$\lambda$-values:
\begin{equation}\label{hintro}
\lambda(ag_1^{4p} + bug_2^{4p}) = 
\begin{cases}
-1 & \text{ if } \deg g_1 \le \deg g_2,\\
1 & \text{ if } \deg g_1 > \deg g_2, 
\end{cases}
\end{equation}
for any $g_1, g_2 \in \kappa[u]$ not both zero,
using the convention $\deg 0 = -\infty$. 
(If $p \ne 2$ then (\ref{hintro}) 
is a special case of Lemma \ref{glemma}, replacing 
$g_1$ and $g_2$ in that lemma with their squares.
We omit the additional
considerations that are
required to verify (\ref{hintro}) 
when $p = 2$.)
In particular, $\lambda(ag_1^{4p} + bug_2^{4p})$ only depends on
the sign of $\ord_\infty(g_1/g_2) = \deg g_2 - \deg g_1$.
The proof of Theorem \ref{mainthm}
rests on a similar counterexample to Chowla's two-variable conjecture, 
with exponent $2p$ rather than $4p$. (See (\ref{eq177}).)
The 
failure of Chowla's conjecture in positive characteristic 
was our initial clue to the possibility that 
elevated rank may occur in non-isotrivial families 
in the function field case.

\vfill

\end{document}